\documentclass[11pt]{amsart}
\usepackage{amsfonts,  amsfonts, amscd, enumerate}

\newtheorem{theorem}{Theorem}[section]

\newtheorem{lemma}[theorem]{Lemma}
\newtheorem{corollary}[theorem]{Corollary}
\newtheorem{proposition}[theorem]{Proposition}
\theoremstyle{definition}
\newtheorem{definition}[theorem]{Definition}

\newtheorem{example}[theorem]{Example}
\newtheorem{examples}[theorem]{Examples}

\theoremstyle{remark}
\newtheorem{remark}[theorem]{Remark}
\newtheorem{remarks}[theorem]{Remarks}
\numberwithin{equation}{section}
\usepackage{amscd}
\usepackage{xypic}
\usepackage{amsmath, amssymb}

\numberwithin{equation}{subsection}
\newcommand{\be}%
  {\protect\setcounter{equation}{\value{subsubsection}}}
  \newcommand{\ee}%
   {\protect\setcounter{subsubsection}{\value{equation}}}

  {\protect\setcounter{subsubsection}{\value{equation}}}

\def \A{{\mathcal A}}

\newcommand{\oC}{\operatorname{Ch}}

\def \C{\mathcal C}

\def \Cl{\mathbb C}

\def \colimm{\underset {m \rightarrow \infty}  {\hbox {lim}}}

\def \colimK.{\underset {\underset K^.  \rightarrow}  {\hbox {lim}}}

\def \colimU.{\underset {\underset U_.  \rightarrow}  {\hbox {lim}}}

\newcommand{\Cqf}{{\rm {Cor}}_{\rm {q.f}}}
\newcommand{\Cf}{{\rm {Cor}}_{\rm {f}}}

\def \D{\mathcal D}

\newcommand{\oDelta}{\operatorname{\Delta}}
\newcommand{\bdD}{\overset {\circ } \Delta }
\newcommand{\nhorn}{{\Lambda _0^n}}

\def \EG1{E{(G \times {\mathbb C}^*)}{\underset {G\times {\mathbb C}^*}  \times}}

\def \EZ(s)1{E{(Z(s) \times {\mathbb C}^*)}{\underset {(Z(s)\times {\mathbb C}^*)}  \times}}

\newcommand{\eps}{ \, {\boldsymbol\varepsilon} \,}

\def \EM(u){EM(u){\underset {M(u)}  \times}}
\def \EM(us){EM(u,s){\underset {M(u, s)}  \times}}

\def \End{{\mathcal E}nd}

\newcommand{\oF}{\operatorname{F}}

\newcommand{\oG}{\operatorname{G}}

\def\holimD{\mathop{\textrm{holim}}\limits_{\Delta }}

\def \Hom{\underline {Hom}}
\def \Hom{{\mathcal H}om}

\def \invlim1{\underset {\infty \leftarrow q}  {\hbox {lim}}^1}

\newcommand{\oLambda}{\operatorname{\Lambda}}
\def \L3{\Lambda \times \Lambda \times \Lambda}
\def \L2{\Lambda \times \Lambda}

\def \longright2arrow{{\overset \longrightarrow  {\overset {}  \longrightarrow}}}

\def \L{L\times \Cl ^*}

\def \M{\mathcal M}

\def \powerl{{\overset l \otimes}}

\def \Q{{\mathbb  Q}}

\def \ra{\rightarrow}

\def \RG^{R(G)^{\hat {}}\ }

\def \res{respectively}

\def \RHom{{{\mathcal R}{\mathcal H}om}}
\def \R{{\mathcal R}}

\def \S{\mathcal S}

\def \topGcoh*{^{top, *} _{G}}
\def \topGho*{ _{top,*} ^{G}}

\def \Z(s){Z(s) \times {\mathbb C}^*}
\def \Z{\mathbb Z}

\begin{document}

\title{Motivic $E_{\infty}$-algebras and the motivic dga}
\author{Roy Joshua}
\address{Department of Mathematics, Ohio State University, Columbus, Ohio, 43210, USA.}
\email{joshua@math.ohio-state.edu}
\thanks{The author was supported by the IHES, the IAS, the MPI and  grants from the NSA and NSF at different stages of this work.}

\begin{abstract}
In this paper we define an explicit ${\rm E}_{\infty}$-structure, i.e. a coherently homotopy
associative and commutative product on chain complexes defining 
(integral and mod-$l$) motivic cohomology as well as mod -$l$ \'etale cohomology. We also discuss several applications.
\end{abstract}

\maketitle

\centerline{\bf Table of contents}
\vskip .3cm 
1. Introduction
\vskip .3cm
2. The basic framework
\vskip .3cm
3. Action of the Barratt-Eccles operads on the motivic complex
\vskip .3cm
4. Mixed Tate motives for a class of smooth schemes over the field
\vskip .3cm
5. Classical cohomology operations
\vskip .3cm
6. Appendix: Chain complexes and operads

\markboth{ROY JOSHUA}{The motivic DGA}
\input xypic

\section{\bf Introduction}
\label{intro}
The main result of this paper is the existence of explicit ${\rm E}_{\infty}$-structures
on complexes defining mod-$l$ motivic and \'etale cohomology of all smooth schemes of finite type over any field. 
Such results have often been implicitly assumed in the past: we make them explicit by constructing
explicit operad actions on these complexes. 
\vskip .3cm
Our constructions already show that there are important
differences between the $E_{\infty}$-structure on the mod$-l$ singular complex of a topological space and the $E_{\infty}$-structure we
provide on the mod-$l$ motivic complex. 
In the first case, this pairing involves the Alexander-Whitney map
from the chain complex associated to the tensor product of two simplicial abelian groups to the tensor product of the
corresponding chain complexes, which is only homotopy commutative in the two arguments. 
In this case, the existence of cohomology operations is a consequence of the fact that the 
Alexander-Whitney map  is only homotopy commutative and not strictly commutative  in the arguments.
In the motivic case, {\it the intersection pairing} on the 
motivic complexes (i.e. the integral motivic complexes and the mod-$l$ motivic complexes) is obtained from an 
external pairing of chain complexes followed by a pairing that involves only shuffle maps which are
 strictly commutative in the arguments of the pairing. (See the beginning of section 3.4 for more details on the
above reasoning.) 
\vskip .3cm
Therefore the source of the lack of strict commutativity for the product on the motivic complex is the lack of strict commutativity of
the corresponding external pairing. 
In other words, our results show that the ${\rm E}_{\infty}$-structure on the singular co-chain complex is a strictly chain-theoretic phenomenon,  whereas the 
${\rm E}_{\infty}$-structure on the motivic complex is not.
The ${\rm E}_{\infty}$-structure on the integral motivic complex also allows the construction of 
 a commutative dga-structure on the rational motivic complex associated to any smooth scheme 
over a field. The latter applies to provide
a construction of categories of mixed Tate motives associated to a large class of schemes.
\vskip .3cm
The following is an outline of the paper. The second section is devoted to setting up the basic framework for the rest of the paper. There we also discuss 
the simplicial
Barratt-Eccles operad. The third section is devoted to defining in detail an action of the simplicial Barratt-Eccles  operad on the integral motivic complexes. 
 This leads to the construction of the motivic dga and 
motivic ${\rm E}_{\infty}$-dga.
We summarize some of these results in the following theorem.
\begin{theorem}
\label{thm2}
  Let $({\rm {Sm}}/k)$ denote the category of all smooth separated schemes of finite type over a fixed field $k$
with morphisms being maps of finite type.  Let $l$ be a prime (not necessarily different from $char (k) =p$) and let $\nu > 0$.
\vskip .3cm
1. Then there exists a functor
\[{\underline \Q}^{mot}: ({\rm {Sm}}/k) ^{op} \ra ( \mbox{DG-algebras over }\Q)\]
\vskip .3cm
\noindent
so that ${\underline \Q}^{mot} = {\underset {r\ge 0} \oplus} {\underline \Q}^{mot}(r)$ and the dga-structure
on ${\underline \Q}^{mot}$ is compatible with the above grading.
If $X$ is any smooth separated scheme of finite type over $k$ and 
$\Q^{mot}_X = R\Gamma (X, {\underline \Q}^{mot})$, then $H^q(\Q^{mot}_X(r)) \cong H_{\mathcal M}^q(X, \Q(r))$ which is the (rational) motivic cohomology of $X$ in degree $q$ and weight $r$.
\vskip .3cm
2. More generally there exist functors
\[{\underline {\mathbb Z}}^{mot}, {\underline {\mathbb Z}}^{mot}/l^{\nu}: ({\rm {Sm}}/k)^{op} \ra (E^{\infty} \mbox{ DG-algebras }), \quad {\underline {\mathbb Z}}^{et}/l^{\nu}:({\rm {Sm}}/k)^{op} \ra (E^{\infty} \mbox{ DG-algebras})\]
\vskip .3cm \noindent
so that 
\[{\underline {\mathbb Z}}^{mot} = {\underset {r\ge 0} \oplus} {\underline {\mathbb Z}}^{mot}(r), \, {\underline {\mathbb Z}}^{mot}/l^{\nu} = 
{\underset {r\ge 0} \oplus} {\underline {\mathbb Z}}^{mot}/l^{\nu}(r) \mbox{ and } {\underline {\mathbb Z}}^{et}/l^{\nu} = {\underset {r\ge 0} \oplus} 
{{\underline {\mathbb Z}}^{et}/l^{\nu}}(r)\]
 with the $E^{\infty}$ dga structures compatible with the above grading. If $X$ is any  smooth separated scheme of finite type over $k$ and 
\[{\mathbb Z}^{mot}_X = R\Gamma (X, {\underline {\mathbb Z}}^{mot}), {\mathbb Z}^{mot}/l_X^{\nu} = R\Gamma (X, {\underline {\mathbb Z}}^{mot}/l^{\nu}) \, 
 \mbox{ (both computed on the Zariski site of X) and}\]
\[{{\mathbb Z}^{et}/l_X^{\nu}} = R\Gamma (X, {{\underline {\mathbb Z}}^{et}/l^{\nu}}) \mbox{ (computed on the \'etale site of X)} \] then,
 \be \begin{align}
H^q({\mathbb Z}^{mot}_X(r)) &\cong H_{\mathcal M}^q(X, {\mathbb Z}(r)) = \mbox{ the motivic cohomology of } X \mbox{ in 
degree q and weight r },\notag \\
H^q({\mathbb Z}^{mot}/l^{\nu}_X(r)) &\cong H_{\mathcal M}^q(X, {\mathbb Z}/l^{\nu}(r)) = \mbox{ the corresponding  mod-l}^{\nu} \mbox {motivic cohomology while } \notag \\
H^q({{\mathbb Z}^{et}/l_X^{\nu}}(r)) &\cong H_{et}^q(X, {\mathbb Z}/l^{\nu}(r)) = \mbox{ the corresponding mod-l}^{\nu} \mbox{ \'etale cohomology}. \notag\end{align}\ee
\vskip .3cm \noindent
 ${\mathbb Z}^{mot}_X$, ${\mathbb Z}^{mot}/l^{\nu}_X$ and ${{\mathbb Z}^{et}/l^{\nu}}_X$ are ${\rm E}_{\infty}$-algebras over the operad
 $N{ Z}(E\Sigma ) \otimes \End_{{\mathcal Z}}$ where $NZ(E\Sigma)$  is the simplicial Barratt-Eccles operad considered in Definition ~\ref{Barratt.Eccles.simp} and $\End_{{\mathcal Z}}$ is the
classical Eilenberg-Zilber operad (discussed in ~\ref{EZ.operad}).
\vskip .3cm
3. Moreover, the functors ${\underline {\mathbb Z}}^{mot}$, ${\underline {\mathbb Z}}^{mot}/l^{\nu}$ and ${\underline {\mathbb Q}}^{mot}$ are
 additive presheaves satisfying cohomological descent (as in  Definition ~\ref{coh.des.MV})
on the big Zariski and Nisnevich sites of smooth quasi-projective schemes over $k$. 
\end{theorem}
\vskip .3cm
Sections four and five of the paper  are devoted to applications. In the fourth section, making use of 
the existence of the motivic dga, we construct a category of mixed Tate motives for a large class of varieties 
that includes all quasi-projective smooth linear varieties over a field which are the complements of closed 
projective smooth linear subvarieties in a bigger projective smooth linear variety.  The main result here is 
Theorem  ~\ref{tatemotives} which we quote presently.
\vskip .3cm
We let $A= {\Q}_X^{mot} = \oplus_{r \ge 0}\Q_X^{mot}(r)$. 
Let ${\mathbb D}_{-}(A)$ denote the derived category of cohomologically bounded below $A$-modules, i.e.
differential graded $A$-modules $M = {\underset {r} \oplus} M(r)$ so that ${\mathcal H}^q(M)(r) =0$ for all sufficiently small $q$. 
 ${\mathcal F}{\mathcal H}_A$ will denote a full sub-category of
${\mathbb D}_{-}(A)$ defined in section  four. ${\mathcal M}{\mathcal T}{\mathcal F}(X)$
is a candidate for the category of mixed Tate motives for the scheme $X$. (See section four for 
terminology and more details):
\begin{theorem}
\label{tatemotives.0}
 If the DGA $A$ is connected (in the sense of ~\ref{hyp.dga.1}), ${\mathcal M}{\mathcal T}{\mathcal F}(X)$ 
is equivalent to the category ${\mathcal F}{\mathcal H}_A$. In particular, this holds for the following 
classes of smooth quasi-projective varieties assuming the Beilinson-Soul\'e conjecture 
holds for the rational motivic cohomology of $Spec \, k$, for example if $k$ is a number field (see Corollary ~\ref{Mixed.Tate.2}):
\vskip .3cm 
(i) all smooth (connected) projective
linear varieties over $k$
\vskip .3cm
(ii) any of the varieties appearing in Example ~\ref{egs} which are also connected, projective and smooth
\vskip .3cm
(iii) any quasi-projective variety $U$ of the form $X-Y$, where $X$ and $Y$ are smooth projective 
varieties as in (i) or (ii) and $Y$ is closed in $X$.
\end{theorem}
\vskip .3cm 
So far the only construction of
a category of mixed Tate motives is for a field and also for the ring of integers in a number field. (See for example, \cite{bl2}, \cite{blk}).
\vskip .3cm
In the fifth section, we show that the operad actions defined earlier lead to {\it classical} cohomology
operations on both mod$-l$ and mod-$p$ motivic cohomology. These operations differ from the operations
constructed by Voevodsky (see \cite{vv}) in the way they behave with respect to the weights. We believe, the methods of this 
paper are the easiest and quickest means of constructing the above classical operations.
 A follow-up to this paper, worked out jointly with
Patrick Brosnan, (See \cite{broj}) explores the relation between the classical cohomology operations as constructed in this paper and the motivic cohomology operations of Voevodsky.
\vskip .3cm
In an earlier draft of the paper, we had provided a different $E_{\infty}$-structure on the complex defining mod$-l$ \'etale cohomology as well as on 
the motives of smooth schemes. This $E_{\infty}$-structure was provided by an algebraic variant of the Eilenberg-Zilber operad. However, this made it 
necessary to first set up and discuss in detail a model structure on the category of chain-complexes of sheaves of $Z/l$-modules on the big Nisnevich or \'etale site
of smooth schemes and seemed to take the focus away from the main results on the action of the simplicial Barratt-Eccles operad on the motivic complexes. Therefore
we have removed all discussions on this second construction which will be discussed elsewhere separately.
\vskip .3cm
It may also be worthwhile pointing out that though the statements corresponding to Theorem ~\ref{thm2} for the singular co-chain complex in algebraic topology 
had been known for a very long time, the explicit construction of an algebra structure for the singular co-chain complex of a topological space over the
simplicial Barratt-Eccles operad is relatively recent: see \cite{JRS} and \cite{bfres}.
It is important for some of our applications, for example the construction of the cohomology
operations in section 5 that are {\it unstable with respect to weight-suspensions}, that we work unstably. Therefore, the $E_{\infty}$-structure
we provide here is an {\it unstable} one on the integral motivic complex,  and distinct from the ring structure on the motivic Eilenberg-Maclane spectrum as a 
${\mathbb P}^1$-spectrum. Moreover,  there are
important differences between this structure and the $E_{\infty}$-ring structure on the singular chain complex of a topological space as 
explained earlier: these justify a careful construction of an explicit $E_{\infty}$-structure on the motivic complex as we do in the present paper.

\vskip .3cm
{\it Acknowledgments}. This has been a rather long project for us, partly because the area of operads had been new to us when we embarked on this 
project and partly because we had been busy with several other more pressing projects.  Over the years, we have benefited from discussions with several mathematicians: we thank Arvind Asok, 
Spencer Bloch, Patrick Brosnan and Zig Fiedorowicz 
 for several helpful discussions/correspondence.  Thanks are also due to an unknown referee for his/her help in improving the exposition.
\section{\bf The basic framework}
\vskip .3cm
Throughout the paper $k$ will denote a fixed field of arbitrary characteristic $p \ge 0$. $({\rm {Sm}}/k)$ will denote the category of separated
 smooth schemes of finite type over $k$. $({\rm {Sm}}/k)_{Zar}$, $({\rm {Sm}}/k)_{Nis}$ and $({\rm {Sm}}/k)_{et}$ will denote this category provided with the 
big Zariski, Nisnevich or \'etale topologies. If $X$ denotes a separated smooth scheme of finite type over $k$, $X_{Zar}$ ($X_{Nis}$, $X_{et}$) will denote 
the big Zariski site (the big Nisnevich site, the big \'etale site, \res) of $X$. we may denote any of these generically by $X_{st}$. 
(An object of $X_{Nis}$ will be a smooth scheme of finite type over $k$,
provided with a map to $X$. Morphisms between two such objects will be commutative triangles over $X$ and coverings will be coverings in the
Nisnevich topology. The sites $X_{et}$ and $X_{Zar}$ may be defined similarly.)
We will let $\S$ denote any one of these sites: since the schemes we consider are of finite type over $k$, it follows that these sites are all skeletally small. $\oC(\S)$ will denote the category of unbounded co-chain complexes of abelian sheaves on $\S$ with differentials of degree $+1$; 
complexes of abelian sheaves with differentials of degree $-1$ will be referred to as
chain complexes. By default, a complex will usually mean a co-chain complex, i.e. one whose differentials are of degree $+1$. The category of all pointed simplicial sheaves on a site $\S$ will be denoted ${\rm {SSh}}(\S)$: observe that
this category contains as a full sub-category, the category of all sheaves of pointed sets. The base point in this
category will be denoted $*$.
\vskip .3cm 
For the most part 
$\oC_{-}(\S)$ ($\oC_0(\S)$) will denote the full sub-category of bounded above complexes
(that are also trivial in positive degrees, \res). For a fixed sheaf of commutative rings $\R$ with unit, $\oC(\S, \R)$ ($\oC_{-}(\S, \R)$, $\oC_0(\S, \R)$) will denote the corresponding categories of complexes of sheaves of $\R$-modules. If $R$ is a fixed commutative ring with unit, we will let $R$ also denote the obvious associated constant sheaf.
Observe that $\oC(\S)$ and $\oC(\S, \R)$  are symmetric monoidal with product $\otimes$ and an internal Hom functor we denote by $\Hom$.  Observe that $\oC_{-}(\S, \R)$ and $\oC_0(\S, \R)$ are not closed under the formation  of the internal hom: this is the main reason for considering the category $\oC(\S, \R)$ of unbounded 
complexes in this paper. While it is convenient  for us to have a model category structure
on $\oC(\S, \R)$ that is compatible with the above tensor structure, such a model structure does not play any key role for the
constructions in the rest of this paper. Therefore, we have chosen not to consider such structures in this paper. 
Throughout the rest of the paper, we will assume that $\R = R$ the constant sheaf associated to a commutative Noetherian ring 
$R$ with unit.
We conclude this section by recalling the definition of the motivic complexes and of the simplicial Barratt-Eccles operad.
\begin{definition}
 Given a fixed $X \eps {\rm {Sm}}/k$, one defines presheaves with transfers,
${\mathbb Z}_{eq}(X)$, ${\mathbb Z}_{tr}(X)$ by $\Gamma (V, {\mathbb Z}_{eq}(X)) = \Cqf(V, X)$ which is
the free abelian group of correspondences on $V \times X$ which are quasi-finite and dominant
over $V$
and $\Gamma (V, {\mathbb Z}_{tr}(X)) = {\rm {Cor}}_{\rm f}(V, X)$ which is the free abelian group of correspondences on $V \times X$ which are finite and surjective over $V$. 
 It is shown in \cite{mvw}, that these define {\it sheaves with transfers} on $({\rm {Sm}}/k)_{Nis}$ and
$({\rm {Sm}}/k)_{Zar}$. 
\end{definition}
\begin{remark}
 Note that ${\mathbb Z}_{eq}(X)$ is more
often denoted $z_{equi}(X, 0)$ in the literature, but the above notation seems more compact and compatible
with ${\mathbb Z}_{tr}(X)$ which is standard in the literature.
\end{remark}

\begin{definition}
\label{Z}
 We let ${\underline {\mathbb Z}}(n) $ denote the complex of sheaves 
$C^*({\mathbb Z}_{tr}(({\mathbb P}^ 1)^{\wedge^n}))[-2n]$ as defined in
~\eqref{motivic.compl.2} on either of the sites $({\rm {Sm}}/k)_{Nis}$ or $({\rm {Sm}}/k)_{Zar}$. 
If $X$ is a smooth scheme over $k$, we also let
${\underline {\mathbb Z}}_X(n)$ denote the restriction of
${\underline {\mathbb Z}}(n)$ to the Zariski or Nisnevich site of $X$.
${\mathbb Z}_X(n)= R\Gamma (X, {\underline {\mathbb Z}}_X(n))$. We 
define ${\underline {\mathbb Z}} = 
{\underset n \oplus}{\underline {\mathbb Z}}(n)$ and ${\mathbb Z}_X = R\Gamma (X, {\underline {\mathbb Z}}_X)$ . ${\underline {\mathbb Z}}/l^{\nu}(n) $ will denote ${\underline {\mathbb Z}}(n) {\underset { Z} \otimes} { Z}/l^{\nu}$ for each fixed prime $l$ and $\nu \ge 0$. ${\underline {\mathbb Z}}^{et}/l^{\nu}(n) $ denotes the same complex on the site $(Sm/k)_{et}$.
\end{definition}
\subsection{\bf The Barratt-Eccles operads }
\vskip .3cm
One of the standard examples of ${\rm E}_{\infty}$-operads are the ones commonly called the Barratt-Eccles operads. We will only consider the simplicial variants of these in this section. (The geometric variant considered in 
an earlier version of the paper will discussed elsewhere as it does not play any role in this paper.)
\vskip .3cm
Given a discrete group $G$, one forms the simplicial group
$EG$ given in degree $n$ by $EG_n = G^{n+1}$. One may verify readily that if $G$ and $H$ are two groups, $E(G \times H)$ is
naturally isomorphic to $EG \times EH$. One may also observe the isomorphism $EG\cong cosk_0(G)$ of simplicial objects: the advantage of $cosk_0(G)$ for us is that it does not involve the group structure. (Recall $cosk_0(G)_m = G^{\times^{m+1}}$. The face map $d_i:cosk_0(G)_m \ra cosk_0(G)_{m-1}$ just drops the $i$-th factor while the degeneracy $s_i:cosk_0(G)_m \ra cosk_0(G)_{m+1}$ just sends the $i$-th factor, $G$, diagonally into $G \times G$ forming  the $i$-th and $i+1$-st factors.) Therefore, taking $cosk_0$, we
see that the map $\gamma_k :\Sigma_k \times \Sigma_{n_1} \times ... \times \Sigma _{n_k} \ra \Sigma _{\Sigma _i n_i}$ defined by $(\sigma _k, \sigma_{n_1}, \cdots, \sigma _{n_k}) \mapsto
 (\sigma_{n_{\sigma_k(1)}}, \cdots, \sigma_{n_{\sigma _k(k)}})$ induces
the structure map $\gamma_k: E(\Sigma _k) \times E(\Sigma _{n_1}) \times E(\Sigma _{n_k}) \ra E(\Sigma _{\Sigma _i n_i})$. 
\vskip .3cm
It may be important to point that since  each $E\Sigma_m$ is identified with $cosk_0(\Sigma_m)$, the operad as defined
above is an operad of simplicial sets only and not of simplicial groups. One may also want to observe that each $E\Sigma_m$ is acyclic and
has a free action by $\Sigma_m$. Given a simplicial set $X$, one lets $X \otimes Spec \, k$ denote
the simplicial scheme defined by $(X\otimes Spec \, k)_n = {\underset {x_n \eps X_n} \sqcup}Spec \, k$ and with the structure maps
induced by the structure maps of the simplicial set $X$. Thus $E\Sigma _m \otimes Spec \, k$ is a simplicial scheme over $Spec \, k$ for each 
$m$. 
\vskip .3cm
\begin{definition} (The simplicial Barratt-Eccles operad)
\label{Barratt.Eccles.simp}
  Henceforth we will let $E\Sigma_m$ denote the simplicial sheaf on the site $({\rm {Sm}}/k)_{Nis}$ represented by this simplicial
scheme. We will apply the free abelian group functor $Z$ followed by the normalization functor $N$ (that produces the Moore-complex: see ~\eqref{N}) to produce an operad
in $\oC(({\rm {Sm}}/k)_{Nis})$ and in $\oC(({\rm {Sm}}/k)_{et})$ from the collection $\{E \Sigma_m |m \ge 0\}$. 
This defines the simplicial Barratt-Eccles operads in $\oC(({\rm {Sm}}/k)_{Nis})$ and $\oC(({\rm {Sm}}/k)_{et})$: this will be denoted $N{  Z}(E\Sigma) = \{N({ Z}(E{\Sigma}_n))|n\}$. Since we do not consider any other variant of this operad in this paper, we will 
refer to this operad henceforth as the {\it Barratt-Eccles operad}. We will often denote this operad by $BE$ with $BE(n) = N({ Z}(E{\Sigma}_n))$.
\vskip .3cm
For later use, it will be important to observe that the pairing induced by $\gamma_k$ on this operad takes the
following form. Given a $q$-simplex, $s_q = (\tau_1, \cdots, \tau_q) \eps N(Z(E \Sigma _k))_q$ and
$p_i$-simplices, $s_{p_i} = (\sigma _1(i), \cdots, \sigma _{p_i}(i)) \eps N(Z(E\Sigma _{j_i}))$, with $j= \Sigma _{i=1}^k j_i$, one first produces the 
$(q, p_1, \cdots, p_k)$-multi-simplex of $N(Z(E\Sigma _j))$ given by
$(\tau_1\circ \sigma _1(1), \cdots, \tau_1 \circ \sigma _{p_1}(1), \cdots , \tau_1 \circ \sigma _1(k), \cdots
\tau_1 \circ \sigma _{p_k}(k), \cdots,  \tau_q\circ \sigma _1(1), \cdots, \tau_q\circ \sigma _{p_1}(1), \cdots , \tau_q \circ \sigma _1(k), \cdots , 
\tau_q \circ \sigma _{p_k}(k))$. (Here we are viewing $\Sigma_k$ acting as the subgroup of $\Sigma_j$ where it permutes the $k$-blocks,
$(1, \cdots , j_1), (j_1+1, \cdots , j_1+j_2), \cdots , (\Sigma_{i=1}^{k-1} j_i +1, \cdots , \Sigma _{i=1}^{k}j_i)$ and each $\Sigma_{j_i}$ is viewed
as the subgroup of $\Sigma _j$ where it permutes the entries of the $i$-th block, $(\Sigma_{l=1}^{i-1}j_l +1, \cdots, \Sigma_{l=1}^{i}j_l)$.) Then one uses the shuffle maps to produce
a $q+ \Sigma _i p_i$-simplex of $N(Z(E \Sigma _j))$. 
\end{definition}
\vskip .3cm
\begin{remark} Recall $M(Spec \, k) = { Z}[0]$ so that $M(E\Sigma _n \otimes Spec \, k) \cong N({Z}(E \Sigma _n))$, where $M$ denotes the functor of {\it motives}
(in the sense of Voevodsky) associated to schemes. \end{remark}
\begin{proposition} The  Barratt-Eccles operad is an acyclic operad in the category $\oC(({\rm {Sm}}/k)_{Nis})$ and $\oC(({\rm {Sm}}/k)_{et})$.  \end{proposition}
\begin{proof}
The acyclicity of the Barratt-Eccles operad follows readily in view of the observation that
each $EG \cong cosk_0(G)$ is acyclic: this simplicial object has an extra degeneracy map $s_{-1}$ induced by the map sending the trivial group $\{e\}$  into $G$ as the identity element of $G$. \end{proof}
\vskip .3cm 
\section{\bf Action of the Barratt-Eccles operads on the motivic complex}
We will assume the basic framework of section 1 throughout this section. In this section, we will define explicitly an action by the simplicial Barratt-Eccles operad 
\newline \noindent
$\{N({ Z}(E\Sigma_n))|n\}$ on the motivic complexes. 
\vskip .3cm
The motivic 
complex ${\underline {\mathbb Z}}(n)[2n]$  will be defined as the complex of quotient sheaves
\be \begin{equation}
\label{motivic.compl.2}
 C^*\Cf( \quad, ({\mathbb P}^1)^{n})/C^*\Cf( \quad, ({\mathbb P}^1)^{n}-{\mathbb A}^n).
\end{equation} \ee
\vskip .3cm
This identification should be clear if we realize that the  complex in ~\eqref{motivic.compl.2}
is nothing but 
\newline \noindent
$C^*{\rm {Cor}}_f(\quad, {\mathbb P}^1 \wedge ... \wedge {\mathbb P}^1)$: recall ${\mathbb P}^1 \simeq {\mathbb G}_m \wedge S^1$ in the 
${\mathbb A}^1$-motivic homotopy category. Therefore, we will denote the complex in ~\eqref{motivic.compl.2} by 
$C^*\Cf (\quad, ({\mathbb P}^1)^{\wedge ^n})$, where ${\mathbb P}^1$ is pointed by $\infty$.
\vskip .3cm
We proceed to define an action of the Barratt-Eccles operad on $\oplus_{n \ge 0}{\mathbb Z}(n)$ with ${\mathbb Z}(n)[2n]$ defined 
as in ~\eqref{motivic.compl.2}. (To gain more insight into the following graph construction
one may consult ~\ref{graph.constr.1}  where we view the motivic complex ${\underline {\mathbb Z}}(n)[2n]$ as 
$C^*({\mathbb Z}_{e.q}(X))$. However, if different definitions of the motivic complex is confusing, the reader should consult
~\ref{graph.constr.1} only after reading through much of the following constructions.)
{\it To motivate and clarify our construction, we will consider first explicitly how one proves (first order) 
homotopy commutativity of the product on the motivic complexes in the context of the above operad action.} 

\subsection{The graph construction.}
\label{graph.constr.3}
There is an obvious action by the group $\Sigma_2$ on 
\[\sqcup_{l, m}({\mathbb P}^1)^{l+m} = \sqcup_{l, m}({\mathbb P}^1)^l \times ({\mathbb P}^1)^m\]
 switching the two factors $({\mathbb P}^1)^l$ and $({\mathbb P}^1)^m$. If $\sigma \eps \Sigma _2$
is the non-identity element, $\sigma ([x_{1,0}: x_{1,1}], \cdots , [x_{l,0}: x_{l, 1}], [y_{1, 0}: y_{1, 1}], \cdots , [y_{m, 0}: y_{m, 1}]) = ([y_{1, 0}: y_{1, 1}], \cdots, [y_{m, 0}: y_{m, 1}], [x_{1,0}: x_{1, 1}], \cdots, 
[x_{l, 0}: x_{l, 1}])$. Let $s_1 =(id, \sigma)$ denote the
obvious $1$-simplex of $E\Sigma _2$.
\vskip .3cm
{\bf The basic strategy} is the following. Let $j_1, j_1'$ and $j_2, j_2'$ denote the integers
$0$ or $1$ so that $j_1+j_1'=1$ and $j_2+j_2'=1$. Observe that 
the ${\mathbb P}^1$ forming the $i_1$-th  factor in $({\mathbb P}^1)^l$  
has homogeneous
coordinates $x_{i_1,0}$ and $x_{i_1,1}$ and the ${\mathbb P}^1$ forming the $i_2$-th  factor in $({\mathbb P}^1)^m$ 
has homogeneous
coordinates $y_{i_2,0}$ and $y_{i_2,1}$. 
 We will begin by defining a map of schemes 
\[{\bar {\phi}}_{s_1}:  {\mathbb A}^1 \times  
(\sqcup _{l, m \ge 1} ({\mathbb P}^1)^{l}  \times ({\mathbb P}^1)^{m}) \ra  \sqcup_{l, m\ge 1}(({\mathbb P}^{1})^{m +l})\]
 by defining its restriction to be the map $\phi_{s_1}$ (as in ~\ref{graph.constr.1}) on the affine space ${\mathbb A}^{l+m} = ({\mathbb A}^1)^l \times ({\mathbb A}^1)^m$ defined by $x_{i_1, j_1'}.y_{i_2, j_2'} \ne 0$ where $i_1$ ranges over $1 \le i_1 \le  l$ and $i_2$ ranges over 
$1 \le i_2 \le m$ with $j_1'$ and $j_2'$ either $0$ or $1$ but depending on $i_1$ and $i_2$, \res. Then we show that these restrictions are compatible under the gluing used to produce
the ${\mathbb P}^1$ forming the various factors in the domain.
 One may in fact define $\bar \phi_{s_1}$ in homogeneous coordinates as follows.
\vskip .3cm
${\bar \phi}_{s_1}(t, [x_{1,0}: x_{1, 1}], \cdots, [x_{l, 0}: x_{l, 1}], [y_{1,0}: y_{1, 1}], \cdots [y_{m, 0}: y_{m, 1}]) =$ the $l+m$-tuple whose $k$-th entry is the point with homogeneous coordinates
\be \begin{align}
&= [tx_{k, j_1}. y_{k, j_2'} +(1-t) y_{k,j2}.x_{k, j_1'}: x_{k, j_1'}.y_{k, j_2'}], k \le min \{l, m\} ,\\
&= [tx_{k, j_1}. x_{k-m, j_2'} +(1-t) x_{k-m,j2}.x_{k, j_1'}: x_{k, j_1'}.x_{k-m, j_2'}], m= min \{l, m\} <k \le l ,  \mbox{ if }m =min \{l, m \}, \notag\\
&=[ty_{k-l, j_1}.x_{k-m, j_2'}+(1-t)x_{k-m, j_2}.y_{k-l, j_1'}: x_{k-m, j_2'}y_{k-l, j_1'}], l <k \le l+m, \mbox{ if }m =min \{l, m \} , \notag\\
&=[ty_{k-l, j_1}. y_{k, j_2'} +(1-t) y_{k,j2}.y_{k-l, j_1'}: y_{k-l, j_1'}.y_{k, j_2'}], l= min \{l, m\} <k \le m, \mbox{ if } l= min \{l, m\}, \notag \mbox{ and} \\
&=[ty_{k-l, j_1}. x_{k-m, j_2'} +(1-t) x_{k-m, j2}.y_{k-l, j_1'}: x_{k-m, j_2'}.y_{k-l, j_1'}], m < k \le l+m, \mbox{ if } l= min \{l, m\} \notag
\end{align} \ee
\vskip .3cm  \noindent
\begin{lemma}
${\bar \phi}_{s_1}$ defines a map ${\mathbb A}^1 \times  (\sqcup _{l, m \ge 1} ({\mathbb P}^1)^{l}  \times ({\mathbb P}^1)^{m}) \ra  \sqcup_{l, m\ge 1}(({\mathbb P}^{1})^{m +l})$  as claimed. 
\end{lemma}
\begin{proof} i.e. The restrictions of the map ${\bar \phi}_{s_1}$ to the affine spaces ${\mathbb A}^1$ which glue together to form the factors ${\mathbb P}^1$ in the {\it domain} are compatible.
We will consider the case $k \le min \{l, m\}$. In this case, it suffices to show for example, the following: let $[x_{k, 0}: x_{k, 1}] \eps {\mathbb P}^1$ forming the $k$-th factor in $({\mathbb P}^1)^l$ and let $[y_{k, 0}: y_{k, 1}] \eps {\mathbb P}^1$ forming the
$k$-th factor in $({\mathbb P}^1)^m$ so that $x_{k,1}. y_{k, 1}\ne 0$ and $x_{k, 0} \ne 0$. Then
the point $[x_{k, 0}: x_{k, 1}]$ is identified with $\frac{x_{k, 0}}{x_{k, 1}}$ in the affine space where $x_{k, 1} \ne 0$ and is identified with $\frac{x_{k, 1}}{x_{k, 0}}$ in the affine space where $x_{k, 0} \ne 0$.
The gluing needed to produce ${\mathbb P}^1$ from these two affine spaces sends 
$$\frac{x_{k, 0}}{x_{k, 1}} \mapsto \frac{x_{k, 1}}{x_{k, 0}}.$$
\vskip .3cm \noindent
We will presently check from the definition that the function ${\bar \phi}_{s_1}$ is defined so as to be compatible with
this identification. There are two affine pieces ${\mathbb A}^1$ covering the ${\mathbb P}^1$ forming the $k$-th factor in $({\mathbb P}^1)^l$: one where $x_{k,0} \ne 0$ and the other where $x_{k, 1} \ne 0$. We are fixing the one of the affine pieces ${\mathbb A}^1$ covering the ${\mathbb P}^1$ forming the $k$-th factor of $({\mathbb P}^1)^m$: we may assume for simplicity that this corresponds to where
$y_{k, 1} \ne 0$. On the affine piece where $x_{k, 1}. y_{k, 1} \ne 0$, 
the $k$-th entry of 
\vskip .3cm
${\bar \phi_{s_1}}(t, [x_{1,0}: x_{1, 1}], \cdots, [x_{l, 0}: x_{l, 1}], [y_{1,0}: y_{1, 1}], \cdots [y_{m, 0}: y_{m, 1}]) = [tx_{k, 0}. y_{k, 1} +(1-t) y_{k,0}.x_{k, 1}: x_{k, 1}.y_{k, 1}]
=t\frac{x_{k,0}}{x_{k, 1}} +(1-t)\frac{y_{k, 0}}{y_{k,1}} $ 
\vskip .3cm \noindent
and on the affine piece where
$x_{k, 0}. y_{k, 1} \ne 0$, the $k$-th entry of 
\vskip .3cm 
${\bar \phi_{s_1}}(t, [x_{1,0}: x_{1, 1}], \cdots, [x_{l, 0}: x_{l, 1}], [y_{1,0}: y_{1, 1}], \cdots [y_{m, 0}: y_{m, 1}]) = 
[tx_{k, 1}. y_{k, 1} +(1-t) y_{k,0}.x_{k, 0}: x_{k, 0}.y_{k, 1}] =
t\frac{x_{k,1}}{x_{k, 0}} +(1-t)\frac{y_{k, 0}}{y_{k,1}}$.
\vskip .3cm \noindent
Therefore, this is compatible with the identification $\frac{x_{k,0}}{x_{k,1}} \mapsto \frac{x_{k,1}}{x_{k,0}}$ used in the gluing to obtain ${\mathbb P}^1$.
There are several remaining cases to consider, which are dealt with similarly.\end{proof}
\vskip .3cm
\begin{example}
To further clarify this definition, one may consider explicitly the  case $l=m=1$.
In this case  ${\bar \phi}_{s_1}: {\mathbb A}^1 \times {\mathbb P}^1 \times {\mathbb P}^1 \ra {\mathbb P}^1 \times {\mathbb P}^1$. Since $l=1$ and $m=1$ we will completely omit the
first subscript of the coordinates. Now the map ${\bar \phi}_{s_1}$ is defined
on the following four affine spaces ($ \cong {\mathbb A}^1 \times {\mathbb A}^1 \times {\mathbb A}^1)$ as follows:
\be \begin{align}
{\bar \phi}_{s_1}(t, [x_0:x_1], [y_0:y_1]) &= ([tx_{1 }. y_{0} +(1-t) y_{1}.x_{0}: x_{0}.y_{0}], [(1-t)x_{1 }. y_{0} +t y_{1}.x_{0}: x_{0}.y_{0}]), \mbox{ if } x_0y_0 \ne 0 ,\\
&= ([tx_{1}. y_1 +(1-t) x_{0}.y_0: x_{0}.y_1], [(1-t)x_{1}. y_1 +tx_{0}.y_0: x_{0}.y_1]), \mbox{ if }x_0y_1 \ne 0, \notag\\
&= ([tx_{0}. y_0 +(1-t) x_{1}.y_1: x_{1}.y_0], [(1-t)x_{0}. y_0 +tx_{1}.y_1: x_{1}.y_0]), \mbox{ if }x_1y_0 \ne 0, \notag \mbox{ and }\\
&=  ([tx_{0}. y_1 +(1-t) x_{1}.y_0: x_{1}.y_1], [(1-t)x_{0}. y_1 +tx_{1}.y_0: x_{1}.y_1]),  \mbox{ if } x_1y_1 \ne 0. \notag 
\end{align} \ee
\vskip .3cm  \noindent
In this case the first (second) factor of ${\mathbb P}^1$ in the domain is obtained by gluing the two affine spaces ${\mathbb A}^1$ corresponding to where $x_0 \ne 0$ and where $x_1\ne 0$
(where $y_0 \ne 0$ and where $y_1 \ne 0$, \res). If $\eta: {\mathbb P}^1_{y_0. y_1 \ne 0} \ra {\mathbb P}^1_{y_1.y_0 \ne 0}$ is the map sending $\frac{y_0}{y_1} \mapsto \frac{y_1}{y_0}$ used
in the gluing to produce ${\mathbb P}^1$, then one verifies that this is compatible with
definition of ${\bar \phi}_{s_1}$ on the first two affine pieces. One verifies the compatibility
of ${\bar \phi}_{s_1}$ with the other gluings used in the domain of ${\bar \phi}_{s_1}$ similarly.
\end{example}
 Let $l, m$ be two fixed positive integers. We begin by
defining a map of presheaves on $({\rm {Sm}}/k)^{op})$
\be \begin{equation}
\label{s.1}
s_1 (=s_1^{l,m}):\Cf(U, ({\mathbb P}^1)^{l+m}) \ra \Cf(U \times \Delta [1], ({\mathbb P}^1)^{l+m})
\end{equation} \ee
\vskip .3cm \noindent
(Strictly speaking the map in ~\eqref{s.1} should be denoted $s_1^{(l,m)}$ as it depends not only on the $1$-simplex $s_1$, but also on the pair $(l, m)$.
But we will usually omit the $(l,m)$ for the sake of brevity.)
This will be contravariant in $U$, so that it is a map of presheaves. Let $Z$ denote a closed
integral subscheme $\subseteq U  \times ({\mathbb P}^1)^{l+m}$ for which the projection $Z \ra U $ is finite.
 Its inverse image under the obvious projection $p: U  \times 
{\mathbb A}^1 \times 
({\mathbb P}^1)^{l+m} \ra U \times ({\mathbb P}^1)^{l+m}$ will define the closed subscheme
$p^{-1}(Z))$. Now the graph of the restriction of $id \times {\bar \phi}_{s_1}$ to this subscheme
defines the closed subscheme $\Gamma _{id \times {\bar \phi}_{s_1}|p^{-1}(Z)}$ of $\Gamma _{id_{U }} \times
\Gamma _{\bar \phi_{s_1}} $ which is contained in $ U \times U  \times {\mathbb A}^1 \times ({\mathbb P}^1)^{l+m} \times ({\mathbb P}^1)^{l+m}$. Since the composition $U  {\overset {\Delta} \ra} U \times U  {\overset {pr_1} \ra} U $ is the identity, and hence in particular proper, the image of 
$\Gamma _{id \times \bar \phi_{s_1}|p^{-1}(Z))}$ under the projection $pr_1 \times id$ to $U \times {\mathbb A}^1 \times ({\mathbb P}^1)^{l+m} \times ({\mathbb P}^1)^{l+m}$ is closed. Since it is the image of an irreducible scheme, it is also irreducible. {\it We will denote this by $s_1'(Z)$.} Since the projection $({\mathbb P}^1)^{l+m} \times ({\mathbb P}^1)^{l+m} \ra ({\mathbb P}^1)^{l+m}$ (which is the projection to the second factor) is clearly proper, we will project $s_1'(Z)$ into $U \times {\mathbb A}^1 \times ({\mathbb P}^1)^{l+m}$ using this map and denote the image
by $s_1(Z)$. We summarize some of the main properties of this construction here:
\vskip .3cm
\subsubsection{}
\label{graph.const.props}
\begin{itemize}
\item
Let $\sigma: ({\mathbb P}^1)^l \times ({\mathbb P}^1)^m \ra ({\mathbb P}^1)^m \times ({\mathbb P}^1)^l$
denote the obvious map interchanging two factors. Then, it is
 clear from the above definition that the scheme $ s_1(Z)$ consists of pairs, consisting of a point $t \eps {\mathbb A}^1$ together with
a point on the  line joining a point of $Z$ with the corresponding point of
$\sigma(Z)$ parametrized by $t \eps {\mathbb A}^1$. Therefore, the projection of
$ s_1(Z))$ to $U \times {\mathbb A}^1$ is in fact {\it finite}. It is also {\it surjective} since the projection of $Z \ra U $ is surjective. 
\item
$p^{-1}(Z)$ being a product of ${\mathbb A}^1$ and $Z$ is evidently  irreducible. Since the graph $\Gamma _{id \times \phi_{s_1}|p^{-1}(Z))}$ is isomorphic to $p^{-1}(Z)$, it is also {\it irreducible}. Therefore, so is $ s_1(Z)$.
\item 
Therefore, $$ s_1(Z) \eps \Cf( U  \times {\mathbb A}^1, ({\mathbb P}^1)^{l+m}).$$
\vskip .3cm \noindent
We will henceforth denote the factor ${\mathbb A}^1$ by $\Delta [1]$. By extending the map $s_1$ by linearity to all cycles in $ \Cf( U  \times {\mathbb A}^1, ({\mathbb P}^1)^{l+m})$, we obtain the required map. 
\item 
The  construction of $s_1(Z)$ is contravariantly functorial in $U$. 
\item
Let $Z_i \eps \Cf(U\times \Delta [n], ({\mathbb P}^1)^{l_i})$, $i=1, \cdots, 2$ be given. We will assume, as before, that each $Z_i$ is, in fact, a closed irreducible subscheme of $U \times \Delta [n] \times ({\mathbb P}^1)^{l_i}$ for which the projection $Z_i \ra U \times \Delta [n]$ is finite.
Let $Z= \Delta ^*(Z_1 \times Z_2)$. Then,  $s_1(\Delta ^*(Z_1 \times Z_2))_{|t=1} =\Delta ^*(Z_1 \times Z_2)$ while,  $s_1(\Delta ^*(Z_1 \times Z_2))_{|t=0} =\Delta ^*(Z_2 \times Z_1)$.
\end{itemize}
\vskip .3cm
Next we proceed to use the map $s_1$ to define a map
\be \begin{equation}
\label{bar.s1}
{\bar s}_1 (={\bar s}_1^{l,m}): \Cf(U, ({\mathbb P}^1)^{\wedge ^{l+m}}) \ra \Cf(U \times \Delta [1], ({\mathbb P}^1)^{\wedge ^{l+m}})
\end{equation} \ee
Starting with the map $s_1$ and passing to the quotient defines a map
$s_1': \Cf(U, ({\mathbb P}^1)^{l+m}) \ra \Cf(U \times \Delta [1], ({\mathbb P}^1)^{\wedge ^{l+m}})$.
We will modify the definition of $s_1'$ so that it will descend to the obvious quotient of the
domain.  We let $\pi_{j}: U \times ({\mathbb P}^1)^{l+m} \ra U \times ({\mathbb P}^1 )^{l+m-1}$ denote 
the projection which omits the $j$-th factor, $j=1, \cdots, l+m$ and let
 $\eta _j:U \times ({\mathbb P}^1)^{l+m} \ra ({\mathbb P}^1)^{l+m} \ra {\mathbb P}^1$ denote 
the composition of the projection to the second factor followed by projection to the $j$-th factor.   
Next let $i_{ j}: ({\mathbb P}^1)^{l+m-1} \ra ({\mathbb P}^1)^{l+m}$, $j =1, \cdots, l=m$ denote the imbedding  with the $j$-th factor
being the base-point $\infty$. (Observe that for any cycle $Z$ in $\Cf(U, ({\mathbb P}^1)^{l+m})$, 
the projection $U \times ({\mathbb P}^1)^{l+m} \ra U$ is finite. Therefore, $\pi_{j}$ defines
a push-forward $\pi_{j*}: \Cf(U, ({\mathbb P}^1)^{l+m}) \ra \Cf(U, ({\mathbb P}^1)^{l+m-1})$.
$i_{j}$ clearly defines a pushforward $\Cf(U, ({\mathbb P}^1)^{l+m-1}) \ra \Cf(U, ({\mathbb P}^1)^{l+m})$.)
Given a closed irreducible and reduced  subscheme $Z \subseteq U \times ({\mathbb P}^1)^{l+m}$, 
we let $D_Z= \{j =1, \cdots, l+m| \eta_j(Z) = \infty \}$. Then we let 
\be \begin{align}
\label{D.Z}
\delta (Z) &=  \Sigma_{j=1}^{l+m} i_{j*} \pi_{j*}(Z), \mbox{ if } D_Z = \phi \\
&= Z+\Sigma _{j  \not \in D_Z} i_{j*} \pi_{j*}(Z), \mbox{ if } D_Z \ne \phi \end{align} \ee
\vskip .3cm \noindent
Now we let 
\be \begin{equation}
{\bar s}_1(Z) = s_1(Z) - s_1(\delta (Z))
\end{equation} \ee
\vskip .3cm
\begin{lemma} Any class in
$\Cf(U, ({\mathbb P}^1)^{l+m} -{\mathbb A}^{l+m})$
 is sent by $\bar s_1$ to zero in the target. The map ${\bar s}_1$ is well-defined and is compatible with
pull-backs in the argument $U$.
\end{lemma}
\begin{proof} 
First let $Z=\infty ^{l+m}$  denote the cycle $U  $ imbedded in $ U \times ({\mathbb P}^1)^{l+m}$ at 
the point $\infty$ in each factor of ${\mathbb P}^1$. Clearly $D_Z = \{1, \cdots, l+m\}$, so that
${\bar s}_1(Z) = s_1(Z) - s_1(Z) = 0$. Next let $Z$ denote a cycle in 
$\Cf(U, ({\mathbb P}^1)^{l+m} -{\mathbb A}^{l+m})$ represented by a closed irreducible and 
reduced  subscheme. Then $D_Z \ne \phi$ and
therefore, ${\bar s}_1(Z) = s_1(Z) - s_1(\delta (Z)) = s_1(Z) - s_1(Z) - 
\Sigma _{j \not \in D_Z} s_1(i_{j*}\pi_{j*}(Z))$. 
\vskip .3cm
For any fixed $j_0 \not \in D_Z$, 
\[s_1(i_{j_{0*}}\pi_{j_{0*}} (Z))= s_1(i_{j_{0*}}\pi_{j_{0*}} (Z)) - s_1(i_{j_{0*}}\pi_{j_{0*}} (Z)) - \Sigma _{j \not \in D_Z, j \ne j_0} s_1(i_{j*} \pi_{j*}i_{j_0*} \pi_{j_0*}(Z)).\] 
For any $j$ appearing
in the last sum, $|\{l| \eta_{l*}(i_{j*}\pi_{j*}i_{j_0*} \pi_{j_0*}(Z)) =\infty \}| > |D_Z|$. Therefore, one may use
ascending induction on $l+m - |D_Z|$ to complete the proof: recall the case $l+m = |D_Z|$ is
when $Z= \infty ^{l+m}$ and this was considered already in the beginning of this proof. 
\vskip .3cm
For a fixed $U$, recall $\Cf(U, ({\mathbb P}^1)^{l+m})$ is the free abelian group generated
by closed irreducible and reduced  subschemes $Z \subseteq U \times ({\mathbb P}^1)^{l+m}$ whose projection
to $U$ is finite and surjective. The well-definedness of the map ${\bar s}_1$ follows from the
observation that ${\bar s}_1$ simply sends the ${\mathbb Z}$-basis elements $Z$ with $D_Z \ne \phi$ to $0$. If $f:U' \ra U$ is a map of smooth schemes of finite type over $k$, any irreducible and reduced  closed subscheme $Z \subseteq U \times ({\mathbb P}^1)^{l+m}$ with $D_Z \ne \phi$ is pulled back to a cycle $f^*(Z) \subseteq U' \times ({\mathbb P}^1)^{l+m}$ with $D_{Z'} \ne \phi$ for each irreducible component $Z'$ of $f^*(Z)$. Therefore, if $Z$ is sent to $0$ by ${\bar s}_1$, then so is $f^*(Z)$. This completes the proof of the lemma.
\end{proof}
\vskip .3cm
Next recall that one has a natural pairing:
\[\Cf(U, ({\mathbb P}^1)^l) \otimes \Cf(U, ({\mathbb P}^1)^m) \ra \Cf(U, ({\mathbb P}^1)^{l+m})\]
\vskip .3cm \noindent
This induces a pairing:
\be \begin{equation}
\label{pairing.11}
\Cf(U, ({\mathbb P}^1)^{\wedge ^l}) \otimes \Cf(U, ({\mathbb P}^1)^{\wedge ^m}) \ra \Cf(U, ({\mathbb P}^1)^{\wedge ^{l+m}}) .
\end{equation} \ee
\vskip .3cm \noindent
Therefore, composing with the map ${\bar s}_1$ considered in ~\eqref{bar.s1}, one obtains a pairing
\be \begin{equation} 
\label{pairing.2}
\mu'({ s_1}, \quad): \Cf(U, ({\mathbb P}^1)^{\wedge ^l}) \otimes \Cf(U, ({\mathbb P}^1)^{\wedge ^m}) \ra \Cf(U \times \Delta [1], ({\mathbb P}^1)^{\wedge ^{l+m}}).
\end{equation} \ee
\vskip .3cm
Observe  the pairing $s_1$ (see ~\eqref{s.1}) is contravariantly functorial in $U$.
Therefore, so is $\bar s_1$ and the  pairing $\mu'({ s_1}, \quad)$. In particular, it follows that
one may replace $U$ by $U \times \Delta [n]$, for any $n$ and that the induced  pairing would
then be compatible with the structure maps of the cosimplicial scheme $\{\Delta [n] |n\}$. 
Therefore,  we obtain a well-defined pairing 
of complexes of quotient presheaves. Clearly this induces a corresponding map of the associated complexes of
quotient sheaves. 
\vskip .3cm
In view of the property, 
\be \begin{equation}
     \label{first.order.homotopy}
s_1(\Delta ^*(Z_1 \times Z_2))_{|t=1} =\Delta ^*(Z_1 \times Z_2) \mbox{ while },  s_1(\Delta ^*(Z_1 \times Z_2))_{|t=0} =\Delta ^*(Z_2 \times Z_1),
\end{equation} \ee
the above construction {\it provides the first order homotopy for
the pairing of the motivic complexes}.
\vskip .3cm 
Observe that, in the description of the first order homotopy as in ~\eqref{first.order.homotopy}, we only considered the
$1$-simplices of $E\Sigma_2$ of the form $(id, \sigma)$, $\sigma \eps \Sigma_2$. 
We may extend the above construction to define an action by all $1$-simplices of $E\Sigma_2$ which are of the form
$(\sigma_1, \sigma_2)$, $\sigma_i \eps \Sigma _2$ as follows: we simply replace the subscheme
$\Delta^*(Z_1 \times Z_2)$ by $\sigma_1^*(\Delta^*(Z_1 \times Z_2))$ where $\sigma_1^*$ is the 
map induced by the permutation action $\sigma_1: \sqcup_{l,m}{\mathbb A}^l \times {\mathbb A}^m \ra \sqcup_{l,m} {\mathbb A}^l \times {\mathbb A}^m$ and apply the same constructions as before
with $\sigma_2 $ playing the role of $\sigma$. We may also define an action of the $0$-simplices
of $E\Sigma_2$ which are given by $\sigma \eps \Sigma_2$ on the motivic complex by sending the
cycle $\Delta^*(Z_1\times Z_2)$ to $\sigma ^*(\Delta^*(Z_1 \times Z_2))$. 
Therefore, the above construction provides a pairing:
\be \begin{multline}
\label{key.pairing.1}
\mu': {\underset s \oplus} Z(s) \otimes \{\Cf(U \times \Delta[n], ({\mathbb P}^1)^{\wedge ^l})|n\} \otimes \{\Cf(U \times \Delta [n], ({\mathbb P}^1)^{\wedge ^m})|n\} \\
\ra \{\Cf(U \times \Delta [n] \times \Delta [1], 
({\mathbb P}^1)^{\wedge ^{l+m}})|n\}
\end{multline} \ee
\vskip .3cm \noindent
where the direct sum is taken over all the $0$- and $1$-simplices in $E\Sigma _2$ and $Z(s)$ 
denotes the chain complex obtained by first taking the (free) simplicial abelian group 
on the $1$-simplex $s$ and then by normalizing it.
\vskip .3cm
This pairing is clearly compatible with restriction to the faces of the $n$-simplex $\Delta [n]$. 
 Moreover, restricting to the two faces of $\Delta [1]$ in $\Delta [n] \times \Delta[1] $ provides the two classes corresponding 
to $\sigma_i^*(\Delta^*(Z_1 \times Z_2))$, $i=0,1 $ as observed in ~\eqref{first.order.homotopy}.  
Observe that both the left and right-sides of the pairing ~\eqref{key.pairing.1} are double-complexes whose bi-degrees are indexed by $n$ and the degree of terms in the complex $Z(s)$. The pairing is compatible with restrictions 
to the faces of $\Delta [n]$ as well as to the faces of the $1$-simplices $ s_1 \eps E\Sigma _2$ as the above descriptions show. Therefore, on taking the associated total complexes, the induced pairing is compatible with the differentials of the complexes on either side. 
The constructions so 
far may now be viewed as a pairing: 
\[\mu':NZ(sk_1E(\Sigma _2)) \otimes {\underline {\mathbb Z}}(l)  \otimes  {\underline {\mathbb Z}}(m) \ra {\underline {\mathbb Z}}(l+ m)\]
\vskip .3cm \noindent
(Here $Z(sk_1E(\Sigma_2))$ denotes the simplicial abelian group obtained by applying the free abelian group functor to the simplicial set $sk_1E(\Sigma_2)$.
 Recall once again that the complex $N(Z(sk_1E(\Sigma_2)))$ is trivial in degrees greater than $1$.)
\subsubsection{Action by higher dimensional simplices}
\vskip .3cm
We will presently extend the above construction to provide  higher order homotopies for the product structure on the motivic 
complexes. {\it Here we consider a $k$-fold product of the motivic complexes and obtain higher order homotopies for the various
resulting products extending the constructions above.} Observe that for $q \ge 1$, a $q$-simplex, $s_{q}$, of $E\Sigma_k$ is given by a sequence $(\sigma_0, \cdots, \sigma _{q})$ with each $\sigma _i \eps \Sigma _k$. 
\subsubsection{}
\label{higher.order}
There is an obvious action by the symmetric group $\Sigma_k$ on $\sqcup_{l_1, \cdots, l_k}{\mathbb A}^{l_1+ \cdots +l_k} = \sqcup_{l_1, \cdots, l_k}{\mathbb A}^{l_1 + \cdots +l_k} $ permuting the $k$-factors ${\mathbb A}^{l_1}, \cdots, {\mathbb A}^{l_k}$. If $s_q =(\sigma_0, \cdots, \sigma _q)$ denotes a $q$-simplex of $E\Sigma _k$,  $\sigma _i$ will also denote the corresponding self-map of ${\mathbb A}^{l_1+ \cdots +l_k}$ henceforth. We define a map 
of schemes $\phi_{s_q}: {\mathbb A}^{q} \times {\mathbb A}^{\Sigma _{i=1}^k l_i}   \ra   {\mathbb A}^{\Sigma _{i=1}^k l_i}  $ by 
\be \begin{equation}
\label{phi.sq}
\phi_{s_q} = t_0\sigma_0 +t_1\sigma_1 + \cdots + t_{q-1}\sigma_{q-1} +(1-t_0 - \cdots -t_{q-1})\sigma _q
\end{equation} \ee
\vskip .3cm \noindent
where $(t_0, \cdots, t_{q-1}) \eps {\mathbb A}^q$. We will identify ${\mathbb A}^q$ with $\Delta [q]$ sending
$(t_0,\cdots, t_{q-1})$ to $(t_0, \cdots, t_{q-1}, (1-t_0 -\cdots -t_{q-1}))$ and view $\phi_{s_q}$ as a map $\Delta [q] \times {\mathbb A}^{\Sigma _i l_i} \ra {\mathbb A}^{\Sigma _i l_i}$.
\vskip .3cm
Next we proceed to adapt the above construction to define a map
\be \begin{equation}
\label{bar.phi.sq}
{\bar \phi}_{s_q} = {\mathbb A}^q \times ({\mathbb P}^1)^{\Sigma _{i=1}^k l_i}   \ra   ({\mathbb P}^1)^{\Sigma _{i=1}^k l_i}  
\end{equation} \ee
as follows. First we let each permutation $\sigma_i$ act on $({\mathbb P}^1)^{\Sigma _{i=1}^k l_i}$ by letting the blocks $({\mathbb P}^1)^{l_1}, \cdots, ({\mathbb P}^1)^{l_k}$ be permuted by $\sigma _i$.
Let $\alpha (\sigma_i(j), p) $ and $\beta(\sigma _i(j), p)$, for $0 \le i \le q$, $1 \le j \le k$ and
$ 1 \le p \le l_{\sigma _i(j)}$ be functions that take values either $0$ or $1$, so that 
$\alpha (\sigma _i(j), p) + \beta (\sigma_i(j), p) =1$. The homogeneous coordinates of a point  in the first-block of $\sigma _i({\mathbb P}^1)^{\Sigma _{i=1}^k l_i} $ are now defined by 
\vskip .3cm
\[([x_{1, \alpha(\sigma_i^{-1}(1), 1)}^{\sigma _i^{-1}(1)}: x_{1, \beta(\sigma _i^{-1}(1), 1)}^{\sigma_i^{-1}(1)} ], \cdots , 
 [x_{l_{\sigma_i ^{-1}(1)}, \alpha(\sigma_i^{-1}(1), l_{\sigma_i ^{-1}(1))}}^{\sigma _i^{-1}(1)}: x_{l_{\sigma_i ^{-1}(1)}, \beta(\sigma _i^{-1}(1), l_{\sigma_i ^{-1}(1))}}^{\sigma_i^{-1}(1)} ])\]
\vskip .3cm \noindent
Therefore, on the affine space defined by $\beta(\sigma_i^{-1}(1), p) \ne 0 $ for all $1 \le p \le l_{\sigma_i ^{-1}(1)}$, $1 \le i \le q$, this point is defined by the coordinates:
\vskip .3cm
\[(\frac{x_{1, \alpha(\sigma_i^{-1}(1), 1)}^{\sigma _i^{-1}(1)}}{ x_{1, \beta(\sigma _i^{-1}(1), 1)}^{\sigma_i^{-1}(1)} }, \cdots , 
 \frac{x_{l_{\sigma_i ^{-1}(1)}, \alpha(\sigma_i^{-1}(1), l_{\sigma_i ^{-1}(1)}  )}^{\sigma _i^{-1}(1)}}{ x_{l_{\sigma_i ^{-1}(1)}, \beta(\sigma _i^{-1}(1),  l_{\sigma_i ^{-1}(1)})}^{\sigma_i^{-1}(1)} })\]
\vskip .3cm \noindent
Therefore, we define ${\bar \phi}_{s_q}$ on the affine piece that is defined by the conditions
$\beta(\sigma_i^{-1}(j), p) \ne 0 $ for all $1 \le p \le l_{\sigma_i ^{-1}(1)}$, $1 \le i \le q$ , $ 1 \le j \le k$ to be given by
\vskip .3cm 
$\Sigma _{i=0}^{q-1}t_i(\frac{x_{1, \alpha(\sigma_i^{-1}(1), 1)}^{\sigma _i^{-1}(1)}}{ x_{1, \beta(\sigma _i^{-1}(1), 1)}^{\sigma_i^{-1}(1)} }, \cdots , 
 \frac{x_{l_{\sigma_i ^{-1}(1)}, \alpha(\sigma_i^{-1}(1),  l_{\sigma_i ^{-1}(1)})}^{\sigma _i^{-1}(1)}}{ x_{l_{\sigma_i ^{-1}(1)}, \beta(\sigma _i^{-1}(1),  l_{\sigma_i ^{-1}(1)})}^{\sigma_i^{-1}(1)} }, \cdots, 
\frac{x_{1, \alpha(\sigma_i^{-1}(k), 1)}^{\sigma _i^{-1}(k)}}{ x_{1, \beta(\sigma _i^{-1}(k), 1)}^{\sigma_i^{-1}(k)} }, \cdots , 
 \frac{x_{l_{\sigma_i ^{-1}(k)}, \alpha(\sigma_i^{-1}(k),  l_{\sigma_i ^{-1}(k)})}^{\sigma _i^{-1}(k)}}{ x_{l_{\sigma_i ^{-1}(k)}, \beta(\sigma _i^{-1}(k),  l_{\sigma_i ^{-1}(k)})}^{\sigma_i^{-1}(k)} })$
\vskip .3cm
$+ (1- \Sigma _{i=0}^{q-1}t_i) (\frac{x_{1, \alpha(\sigma_q^{-1}(1), 1)}^{\sigma _q^{-1}(k)}}{ x_{1, \beta(\sigma _q^{-1}(k), 1)}^{\sigma_q^{-1}(k)} }, \cdots , 
 \frac{x_{l_{\sigma_q ^{-1}(k)}, \alpha(\sigma_q^{-1}(k),  l_{\sigma_q ^{-1}(k)})}^{\sigma _q^{-1}(k)}}{ x_{l_{\sigma_q ^{-1}(k)}, \beta(\sigma _q^{-1}(k),  l_{\sigma_q ^{-1}(k)})}^{\sigma_q^{-1}(k)} }, \frac{x_{1, \alpha(\sigma_q^{-1}(k), 1)}^{\sigma _q^{-1}(k)}}{ x_{1, \beta(\sigma _q^{-1}(k), 1)}^{\sigma_q^{-1}(k)} }, \cdots , 
 \frac{x_{l_{\sigma_q ^{-1}(k)}, \alpha(\sigma_q^{-1}(k),  l_{\sigma_q ^{-1}(k)})}^{\sigma _q^{-1}(k)}}{ x_{l_{\sigma_q ^{-1}(k)}, \beta(\sigma _q^{-1}(k),  l_{\sigma_q ^{-1}(k)})}^{\sigma_q^{-1}(k)} })$
\vskip .3cm \noindent
It suffices to show that this definition is compatible with the identifications used when one glues together two
copies of ${\mathbb A}^1$ to produce ${\mathbb P}^1$. This may be verified as in the case $k=2$ and $q=1$ that
we considered above.
\vskip .3cm
 For any scheme $U$, the graph of the map, $\Gamma _{id \times {\bar \phi}_{s_q}}$,
is a closed subscheme of $U \times \Delta [q] \times U \times \Delta [q] \times 
({\mathbb P}^1)^{\Sigma _i l_i} \times ({\mathbb P}^1)^{\Sigma _i l_i}$ naturally isomorphic to 
$U \times \Delta [q] \times ({\mathbb P}^1)^{\Sigma _i l_i}$. Let 
$Z \subseteq U \times \Pi_{i=1}^k ({\mathbb P}^1)^{l_i}$ be a closed irreducible and 
reduced sub-scheme so that the projection to $U$ is finite and surjective. (For example, one may begin
 with  $Z_i \eps \Cf(U\times \Delta [n], ({\mathbb P}^1)^{l_i})$, $i=1, \cdots, k$. We will assume, as before, that each $Z_i$ is, in fact, a closed
 irreducible and reduced  subscheme of $U \times \Delta [n] \times ({\mathbb P}^1)^{l_i}$ for which the projection $Z_i \ra U \times \Delta [n]$ is finite.
Let $Z$ denote an irreducible component of  $\Delta ^*(\Pi_{i=1}^k Z_i)$. )
  Its inverse image under the obvious projection $p: U \times \Delta [n] \times \Delta [q] \times
({\mathbb P}^1)^{\Sigma _{i=1}^k l_i} \ra U \times \Delta [n] \times ({\mathbb P}^1)^{\Sigma _{i=1}^k l_i}$ will define the closed subscheme
$p^{-1}(\Delta^*( \Pi_{i=1}^k Z_i))$. Now the graph of the restriction of $id \times {\bar \phi}_{s_q}$ to this subscheme defines the closed subscheme $\Gamma _{id \times {\bar \phi}_{s_q}|p^{-1}(Z ))}$ of $\Gamma _{id_{U \times \Delta [n ]}} \times
\Gamma _{{\bar \phi}_{s_q}} $ which is contained in $ U \times \Delta [n] \times U \times \Delta [n] \times \Delta [q] \times ({\mathbb P}^1)^{\Sigma _{i=1}^k l_i} \times ({\mathbb P}^1)^{\Sigma _{i=1}^k l_i}$. 
\vskip .3cm
Since the composition $U \times \Delta [n] {\overset {\Delta} \ra} U \times \Delta [n] \times U \times \Delta [n] {\overset {pr_1} \ra} U \times \Delta [n]$ is the identity, and hence in particular proper, the image of 
$\Gamma _{id \times {\bar \phi}_{s_1}|p^{-1}(Z )}$ under the projection
 $pr_1 \times id$ to $U \times \Delta [n] \times \Delta [q] \times ({\mathbb P}^1)^{\Sigma_{i=1}^kl_i} \times ({\mathbb P}^1)^{\Sigma _{i=1}^kl_i}$ is closed. Since it is the image of an irreducible scheme, it is also irreducible. {\it We will denote this by 
$s_q'(Z) = s_q'(Z))$.} We will next project this into  $U \times \Delta [n] \times \Delta [q] \times ({\mathbb P}^1)^{\Sigma_{i=1}^kl_i}$ using the  projection 
 $({\mathbb P}^1)^{\Sigma _{i=1}^kl_i}\times ({\mathbb P}^1)^{\Sigma _{i=1}^kl_i} \ra ({\mathbb P}^1)^{\Sigma _{i=1}^kl_i}$
 to the second factor and denote the image of $s_q'(Z) $ by $s_q(Z)$.
 It is clear from the above definition that the projection of
$s_q'(Z)$  to $U  \times \Delta [q]$ is in fact {\it finite}.
It is also {\it surjective} since the projection of $Z \ra U $ is 
 surjective. Therefore, we obtain the map:
\be \begin{equation}
\label{sq.0}
s_q (=s_q^{l_1, \cdots, l_k}): {\underset {s_q \eps sk_qE(\Sigma _k)} \oplus} Z(s_q)_q \otimes \Cf(U, ({\mathbb P}^1)^{\Sigma _{i=1}^k l_i})  \ra  \Cf(U \times \Delta[q], ({\mathbb P}^1)^{\Sigma _{i=1}^kl_i}) \notag
\end{equation} \ee
\vskip .3cm \noindent
where  $Z(s_{q})$ denotes the sub-complex of $ZNE\Sigma _k$ generated by the $q$-simplex $s_{q}$.  One may observe that last pairing defines in fact a map of complexes:
\be \begin{equation}
\label{sq}
s_q: {\underset {s_q \eps sk_qE(\Sigma _k)} \oplus} Z(s_q) \otimes \Cf(U, ({\mathbb P}^1)^{\Sigma _{i=1}^k l_i})  \ra \Cf(U \times \Delta[q], ({\mathbb P}^1)^{\Sigma _{i=1}^kl_i}) \notag
\end{equation} \ee
\vskip .3cm \noindent
As in ~\eqref{bar.s1}, we may use this to define a map
\be \begin{equation}
\label{sq}
{\bar s}_q (={\bar s}_q^{l_1, \cdots, l_k}): {\underset {s_q \eps sk_qE(\Sigma _k)} \oplus} Z(s_q) \otimes \Cf(U, ({\mathbb P}^1)^{\wedge ^{\Sigma _{i=1}^k l_i}})  \ra \Cf(U \times \Delta[q], ({\mathbb P}^1)^{ \wedge ^{\Sigma _{i=1}^kl_i}}) \notag
\end{equation} \ee
\vskip .3cm
Composing with the pairing $\Cf (U, ({\mathbb P}^1)^{\wedge ^{l_1}}) \otimes \cdots \Cf(U, ({\mathbb P}^1)^{\wedge ^{l_k}}) \ra \Cf(U, ({\mathbb P}^1)^{\wedge ^{\Sigma _{i=1}^k l_i}})$, this provides a pairing 
\be \begin{equation}
\label{sq.1}
{\bar s}_q: {\underset {s_q \eps sk_qE(\Sigma _k)} \oplus} Z(s_q) \otimes \otimes_{i=1}^k \Cf(U, ({\mathbb P}^1)^{ \wedge ^{l_i}})  \ra \Cf(U \times \Delta[q], ({\mathbb P}^1)^{\wedge ^{\Sigma _{i=1}^kl_i}})
\end{equation} \ee
\vskip .3cm 
As before, one may verify this pairing is contravariantly functorial in $U$, and is 
compatible with the face maps sending the simplex $s_q$ to its faces, so that one obtains the pairing of complexes:
\be \begin{multline}
\label{key.pairing.2}
\mu'({s_q}, \quad ):{\underset {s_q \eps sk_qE(\Sigma _k)} \oplus} Z(s_q) \otimes \{\Cf(U \times \Delta[n], ({\mathbb P}^1)^{l_1})|n\} \otimes \cdots \otimes \{\Cf(U \times \Delta [n], ({\mathbb P}^1)^{l_k})|n\} \\
\ra \{\Cf(U \times \Delta [n] \times \Delta [q],  ({\mathbb P}^1)^{\Sigma _{i=1}^k l_k})|n\}
\end{multline} \ee
\vskip .3cm
Making use of the isomorphisms defined using shuffle maps as in ~\eqref{shuffle.iso} and taking the sum over all such shuffles, one may similarly define a pairing
\be \begin{multline}
\label{key.pairing.3}
\mu'({s_q}, \quad ): Z(s_{q}) \otimes \Cf(U \times \Delta[n] \times \Delta [p_1], ({\mathbb P}^1)^{\wedge ^{l_1}}) \otimes \cdots \otimes \Cf(U \times \Delta [n] \times \Delta [p_k], ({\mathbb P}^1)^{\wedge ^{l_k}}) \\
\ra \Cf(U \times \Delta [n+ \Sigma _i p_i] \times \Delta [q], ({\mathbb P}^1)^{{\wedge ^{\Sigma _{i=1}^k l_i}}})
\end{multline} \ee
where $Z(s_{q})$ denotes the sub-complex of $ZNE\Sigma _k$ generated by the $q$-simplex
$s_{q}$.  
\vskip .3cm
\subsection{}
\label{symm.act}
Clearly the action of the symmetric group $\Sigma _k$ on the simplex $s_{q}$, $(\sigma', \sigma _1, \cdots \sigma _q) \mapsto   (\sigma _1 \circ \sigma', \cdots,  \circ \sigma _q \circ \sigma')$ corresponds to the action of $\sigma ' \eps \Sigma _k$ on $\Cf(U \times \Delta[n], ({\mathbb P}^1)^{\Sigma _i l_i})$ permuting the weight-factors $({\mathbb P}^1)^{l_1}, \cdots, ({\mathbb P}^1)^{l_k}$. 
\vskip .3cm
Though the above constructions provide higher order homotopies, {\it we need to extend these to 
multi-simplices so as to obtain the associativity of the operad action as proved in}
Theorem ~\ref{assoc}. We begin by extending the above construction to bi-simplices.
\vskip .3cm
If $s_q=(\sigma_0, \cdots, \sigma _q)$, $s_p=(\tau_0, \cdots, \tau_p)$ are $q$ and 
$p$ simplices of $E\Sigma _k$, let   $\phi_{(s_q, s_p)}: \Delta [q] \times \Delta [p] \times 
{\mathbb A}^{\Sigma _{i=1}^k l_i} \ra {\mathbb A}^{\Sigma _{i=1}^k l_i}$ be the map defined by  
\vskip .3cm
$(t_0\sigma_0 +t_1\sigma_1 + \cdots + t_{q-1}\sigma_{q-1} +(1-t_0 - \cdots -t_{q-1})\sigma _q) \circ (s_0\tau_0 +s_1\tau_1 + \cdots + s_{p-1}\tau_{p-1} +(1-s_0 - \cdots -s_{p-1})\tau _p) $
\vskip .3cm
$= \Sigma _{i=0, j=0}^{i=q-1, j=p-1} t_i.s_j \sigma _i \circ \tau _j + (1-t_0 - \cdots -t_{q-1}).\Sigma _{j=0}^{p-1}s_j\sigma _q \circ \tau_j + (1-t_0 - \cdots -t_{q-1}). (1-s_0 - \cdots s_{p-1})\sigma _q \circ \tau_p$  with $ (t_0, \cdots, t_q) \eps {\mathbb A}^q$, $(s_0, \cdots, s_p) \eps {\mathbb A}^p$.
\vskip .3cm \noindent
One then obtains a corresponding map $\bar \phi_{(s_q, s_p)}: \Delta [q] \times \Delta [p] \times 
({\mathbb P}^1)^{\Sigma _{i=1}^k l_i} \ra ({\mathbb P}^1)^{\Sigma _{i=1}^k l_i}$. Applying the graph constructions as above, one obtains a pairing:
\be \begin{multline}
\label{key.pairing.4}
\mu'({s_q, s_p}, \quad ): (Z(s_{q}) \otimes Z(s_p)) \otimes \Cf(U \times \Delta[n] \times \Delta [p_1], ({\mathbb P}^1)^{\wedge ^{l_1}}) \otimes \cdots \otimes \Cf(U \times \Delta [n] \times \Delta [p_k], ({\mathbb P}^1)^{\wedge ^{l_k}}) \\
\ra \Cf(U \times \Delta [n+ \Sigma _i p_i] \times \Delta [q] \times \Delta [p], ({\mathbb P}^1)^{\wedge ^{\Sigma _{i=1}^k l_i}})
\end{multline} \ee
\vskip .3cm \noindent
Since $\phi_{(s_q, s_p)}: \Delta [q] \times \Delta [p] \times {\mathbb A}^{\Sigma _{i=1}^k l_i} \ra 
{\mathbb A}^{\Sigma _{i=1}^k l_i}$ is the composition $\phi_{(s_q)} \circ (id_{\Delta [q]} \times \phi_{(s_p)})$, one may see readily that 
$\bar \phi_{(s_q, s_p)} = \bar \phi_{(s_q)} \circ (id_{\Delta [q]} \times \bar \phi_{(s_p)})$. 
 Observe that the action $ \sigma_i$ and $\tau_j$ on 
${\mathbb A}^{\Sigma _{j=1}^k l_j}$ is linear and hence $\tau_j, \sigma _i$ commute with the variables $t_0, \cdots, t_{q-1}, s_0, \cdots, s_{p-1}$. Moreover, for each increasing map $(\phi, \psi):[p+q] \ra [p] \times [q]$, one may see readily from the  above definition that 
\be \begin{equation}
\label{shuffle.real}
\mu'((\phi, \psi)^*(s_q \times s_p), \quad ) = (\phi, \psi)^* \circ  \mu'(s_q, s_p, \quad ) \circ \eta
\end{equation} \ee
Here $(\phi, \psi)^*(s_q \times s_p) \eps E\Sigma _k$ denotes the $q+p$-simplex defined as in 
(~\ref{shuffle.product}). ($(\phi, \psi)^*:  \Cf(U \times \Delta [n+ \Sigma _i p_i] \times \Delta [q] \times \Delta [p], ({\mathbb P}^1)^{\wedge ^{\Sigma _i l_i}}) \ra  \Cf(U \times \Delta [n+ \Sigma _i p_i] \times \Delta [q+p], ({\mathbb P}^1)^{\wedge ^{\Sigma _i l_i}})$ is defined using the 
shuffle map $(\phi, \psi)$. $\eta: Z((\phi, \psi)^*(s_q \times s_p))  \ra Z(s_q) 
\otimes Z(s_p)$ is the obvious map induced by the isomorphism $Z((\phi, \psi)^*(s_q \times s_p))_{p+q} {\overset {\cong} \ra} Z(s_q)_q \otimes Z(s_p)_p$ and the fact that $Z((\phi, \psi)^*(s_q \times s_p))$ is generated
by the single simplex in degree $p+q$.
\vskip .3cm
One may extend the  pairing ~\eqref{key.pairing.4} to products of several simplices: i.e. if 
$s_{p_1}, \cdots, s_{p_l}$ are
simplices of $E\Sigma _k$ of dimensions $p_1, \cdots, p_l$ and $s_q \eps E\Sigma_k$ is a 
$q$-simplex, then one may define a pairing:
\be \begin{multline}
\label{assoc.key.2}
\mu'(s_q, s_{p_1}, \cdots, s_{p_l}, \quad ): (Z(s_{q}) \otimes Z(s_{p_1})) \otimes \cdots \otimes Z(s_{p_l}))
 \otimes \Cf(U \times \Delta[n], ({\mathbb P}^1)^{\wedge ^{l_1}}) \otimes \cdots \otimes \Cf(U \times \Delta [n], ({\mathbb P}^1)^{\wedge ^{l_k}}) \\
\ra \Cf(U \times \Delta [n] \times \Delta [q] \times \Delta [p_1] \times \cdots \Delta [p_l], ({\mathbb P}^1)^{\wedge ^{\Sigma _{i=1}^k l_i}})
\end{multline} \ee
\vskip .3cm \noindent
Moreover if $\alpha: [q +p_1 + \cdots + p_l] \ra [q] \times [p_1] \times \cdots \times [p_l]$ is any increasing map, then composing the last pairing with $\alpha ^*:\Cf(U \times \Delta [n] \times \Delta [q] \times \Delta [p_1] \times \cdots \Delta [p_l], ({\mathbb P}^1)^{\wedge ^{\Sigma _i l_i}}) \ra
\Cf(U \times \Delta [n] \times \Delta [q + p_1 + \cdots  + p_l], ({\mathbb P}^1)^{\wedge ^{\Sigma _{i=1}^k l_i}})$ and pre-composing this with the (obvious) map
 $Z(\alpha ^*(s_q \times s_{p_1} \times \cdots \times s_{p_l})) \ra Z(s_q) \otimes Z(s_{p_1}) 
\otimes \cdots \otimes Z(s_{p_l})$ identifies with the  pairing $\mu'({\alpha ^*(s_q, s_{p_1}, \cdots, s_{p_l}}, \quad ))$.
\vskip .3cm
\subsubsection{}
\label{pairing.1.shuffle}
One may observe that the  construction of the higher order homotopies as in ~\eqref{assoc.key.2} takes place in the weight-factor 
$({\mathbb P}^1)^{\Sigma _jl_j}$ only and is (contravariantly) functorial in the remaining arguments: it follows readily that the  pairing in ~\eqref{assoc.key.2}  is compatible with pull-backs in the remaining arguments. This, in turn, implies that the above pairing commutes with shuffle maps (which are induced by pull-backs on the non-weight-factors) 
in the following sense so that the composition (where $s_q$ denotes a $q$-simplex of $NZ(E\Sigma _k)$)
\vskip .3cm 
$Z(s_q) \otimes \Cf(U \times  \Delta [n] \times \Delta [p_1], ({\mathbb P}^1)^{\wedge ^{l_1}}) \otimes \cdots \otimes \Cf(U \times \Delta [n] \times \Delta [p_k], ({\mathbb P}^1)^{\wedge ^{l_k}}) 
 {\overset {\mu'} \ra} \Cf(U \times \Delta [n] \times \Delta [p_1]\times \cdots \times \Delta [p_k] \times \Delta [q], ({\mathbb P}^1)^{\wedge ^{\Sigma _{i=1}^k l_i}}) {\overset {{\rm {shuffle}}} \ra } 
\Cf(U \times \Delta [n+\Sigma _ip_i] \times \Delta [q], ({\mathbb P}^1)^{\wedge ^{\Sigma _{i=1}^k l_i}})$
\vskip .3cm \noindent
factors as
\vskip .3cm
$Z(s_q) \otimes \Cf(U \times  \Delta [n] \times \Delta [p_1], ({\mathbb P}^1)^{\wedge ^{l_1}}
) \otimes \cdots \otimes \Cf(U \times \Delta [n] \times \Delta [p_k], ({\mathbb P}^1)^{\wedge ^{l_k}})
{\overset {id \otimes {\rm {shuffle}}} \longrightarrow } Z(s_q) \otimes  \Cf(U \times  \Delta [n+\Sigma _ip_i] , ({\mathbb P}^1)^{\wedge ^{l_1}})
 \otimes \cdots \otimes \Cf(U \times \Delta [n+\Sigma _ip_i] , ({\mathbb P}^1)^{\wedge ^{l_k}})
{\overset {\Delta ^* \circ \mu'} \ra } \Cf(U \times \Delta [n+\Sigma _ip_i] \times \Delta [q], ({\mathbb P}^1)^{\wedge^{\Sigma _{i=1}^k l_i}})$.
\vskip .4cm \noindent
Composing the above pairing with another (obvious) shuffle map, one maps the last term above to 
$\Cf(U \times \Delta[n+\Sigma _i p_i +q], ({\mathbb P}^1)^{\wedge ^{\Sigma _{i=1}^k  l_i}})$.
\vskip .3cm 
We may also start with cycles $Z_i \eps \Cf(U \times \Delta[n_i], ({\mathbb P}^1)^{\wedge ^{
l_i}})$, $i=1, \cdots, k$; using shuffle maps
one first produces corresponding cycles in 
$\Cf(U \times \Delta[n], ({\mathbb P}^1)^{\wedge ^{\Sigma _{i=1}^kl_i}}
)$ where  $n = \Sigma _i n_i$ and then one applies the above construction. Moreover, the
above description shows the action is compatible with restriction to the
faces of the simplex $s_{q}$, i.e. the $\Delta [q]$ appearing above. It is also compatible with restriction to the
faces of $\Delta [n]$ and $\Delta [p_i]$ appearing above, so that the last pairing defines a pairing of complexes:
\be \begin{equation}
\label{pairing.2}
NZ(s_{q}) \otimes {\underline {\mathbb Z}}(l_1) \otimes \cdots \otimes {\underline {\mathbb Z}}(l_k) \ra {\underline {\mathbb Z}}(\Sigma _{i=1}^kl_i)
\end{equation} \ee
\vskip .3cm \noindent
where ${\underline {\mathbb Z}}(m)$ {\it here }denotes the motivic complex of weight $m$. $Z(s_{q})$ denotes the sub-chain complex of 
$NZ(E\Sigma _k)$ generated by a {\it chosen} $q$ simplex $s_{q}$. 
Let $\Delta [n]_{ss}$ denote the usual simplicial set given by $\Delta [n]_{ss,k} =Hom_{\Delta}([k], [n])$. By composing with the obvious map 
$\Delta [s_{q}] \ra E\Sigma _k$ sending the
generator $i_{s_{q}}$ to the simplex $s_{q}$, we also obtain a pairing
\be \begin{equation}
\label{pairing.3}
NZ(\Delta[q]_{ss})) \otimes {\mathbb Z}(l_1) \otimes \cdots \otimes {\mathbb Z}(l_k) \ra {\mathbb Z}(\Sigma _{i=1}^kl_i)
\end{equation} \ee
\vskip .3cm \noindent
where $NZ(\Delta[q]_{ss})$ is the obvious chain complex obtained from the simplicial abelian group 
$Z(\Delta [q]_{ss})$.
\vskip .3cm
By ascending induction on $q \ge 1$, we may now assume that we have already established a pairing
\[NZ(sk_{q-1}(E\Sigma _k)) \otimes {\underline {\mathbb Z}} \otimes \cdots \otimes {\underline {\mathbb Z}} \ra {\underline {\mathbb Z}}, \mbox{ where }
{\underline {\mathbb Z}} = \oplus _r {\underline {\mathbb Z}}(r)\]
\vskip .3cm \noindent
denotes the motivic complex and there are $k$-factors of this on the left. $NZsk_{q-1}(E\Sigma _k)$ denotes the chain complex obtained by normalizing the simplicial abelian group obtained by
applying the free abelian group functor $Z$ dimension-wise to $sk_{q-1}E\Sigma _k$. 
Observe the $q$-skeleton of $E\Sigma _k$ is the filtered colimit of all its $m$-cells, $ m \le q$. Now observe the co-cartesian square (where the sum ranges over all $q$-simplices of $NZ(E\Sigma _k)$): 
\[
\xymatrix{{\oplus {NZ({\overset {\circ }\Delta}[q]_{ss})}} \ar@<1ex>[r] \ar@<-1ex>[d] & {\oplus {NZ(\Delta [q]_{ss})}} \ar@<1ex>[d]\\
 {NZ(sk_{q-1}  (E\Sigma _k))} \ar@<1ex>[r] & {NZ(sk_q(E\Sigma _k))}}\]
\vskip .3cm \noindent
  Then one uses ascending induction on $q$, the above co-cartesian square  and the pairing ~\eqref{pairing.3} to define the  pairings 
\be \begin{equation}
\label{pairing.4}
NZ(sk_q(E\Sigma _k)) \otimes {\underline {\mathbb Z}}(l_1) \otimes \cdots \otimes  
{\underline {\mathbb Z}}(l_k) \ra {\underline {\mathbb Z}}(\Sigma _{i=1}^k l_i)
\end{equation} \ee
where there are $k$-factors of ${\mathbb Z}$ on the left and which are compatible with respect to the skeletal filtration on $E\Sigma_k$.  
(Clearly one may start the induction when $q=1$.) Finally take the colimit over $q \ra \infty$ to obtain the pairing
\be \begin{equation}
\label{pairing.5}
\mu_k:NZ(E\Sigma _k) \otimes {\underline {\mathbb Z}}(l_1) \otimes \cdots \otimes  
{\underline {\mathbb Z}}(l_k) \ra {\underline {\mathbb Z}}(\Sigma _{i=1}^k l_i)
\end{equation} \ee
\subsection{}
\label{mu.comm}
The observation in ~\ref{symm.act} shows that the action of the symmetric
group $\Sigma _k$ by an element $\sigma \eps \Sigma _k$ on $NZ(E\Sigma _k)$
cancels out with the action by the element $\sigma ^{-1}$ on the
$k$-factors ${\mathbb Z}(l_1), \cdots, {\mathbb Z}(l_k)$. 
Abbreviating ${\underset {l} \oplus} {\underline {\mathbb Z}}(l)$ to ${\underline {\mathbb Z}}$, the above pairing may be shortened to $\mu_k:NZ(E(\Sigma _k)) \otimes {\underline {\mathbb Z}}^{\otimes ^k} \ra {\underline {\mathbb Z}}$.
\vskip .3cm
\subsubsection{}
\label{bullet.not}
Given cycles $Z_i \eps \Cf(U \times \Delta [n], ({\mathbb P}^1)^{\wedge ^{l_i}})$, $i=1, \cdots, s$, we let
\[{\underset {i=1} {\overset s \bullet}} Z_i = \Delta ^*( Z_1 \times \cdots Z_s) \eps \Cf(U \times \Delta [n], ({\mathbb P}^1)^{\wedge^{\Sigma _i l_i}}). \]
\vskip .3cm 
Next one starts with cycles $Z(j_i)_1, \cdots Z(j_i)_{j_i}$ in $\Cf(U \times \Delta [ n], ({\mathbb P}^1)^{\wedge^{l(j_i)_1}}), \cdots, \Cf (U \times \Delta [n], ({\mathbb P}^1)^{\wedge^{l(j_i)_{j_i}}})$, $i=1, \cdots, k$. Let $j= \Sigma _{i=}^k j_i$ and $s_{j_1},\cdots, s_{j_k}$ denote simplices of dimension $p_1, 
\cdots, p_k$ in $E\Sigma_{j_1}, \cdots, E\Sigma_{j_k}$, \res. Now one may observe readily that
\be \begin{align}
\label{pairing.5}
\mu'(s_{j_1}, \cdots,  s_{j_k}, {\underset {i=1} { \overset {j_1} \bullet}}Z(j_1)_i,  \cdots 
, {\underset {i=1} { \overset {j_k} \bullet}}Z(j_k)_i) &= \mu'(s_{j_1}, 
({\underset {i=1} { \overset {j_1} \bullet}}Z(j_1)_i)) \bullet \cdots \bullet \mu'(s_{j_k}, 
({\underset {i=1} { \overset {j_1} \bullet}}Z(j_k)_i)\\
&= \mu'(s_{j_1}, Z(j_1)_1, \cdots, Z(j_1)_{j_1}) \bullet \cdots \bullet \mu'(s_{j_k}, Z(j_k)_1, \cdots, 
Z(j_k)_{j_k})
\end{align} \ee
\vskip .3cm \noindent
The left-hand-side is defined as in ~\eqref{assoc.key.2} by using the obvious diagonal imbedding  
$\Sigma _{j_1} \times \cdots \times \Sigma_{j_k}$ in $\Sigma _j$ so that all the 
simplices $s_{j_i}$ are viewed as simplices in $E\Sigma_j$. In view of this, the actions of
the simplices $s_{j_i}$ commute with each other, with $s_{j_i}$ acting trivially
on all the terms ${\underset {i=1} { \overset {j_1} \bullet}}Z(j_l)_i, l \ne i$, 
so that one obtains the identifications stated
above. The right-hand-side defines 
the image of the product
$ \mu'(s_{j_1}, {\underset {i=1} { \overset {j_1} \bullet}}Z(j_1)_i) \bullet \cdots \bullet 
\mu'(s_{j_k}, {\underset {i=1} { \overset {j_k} \bullet}}Z(j_k)_i)$ in 
$\Cf (U \times \Delta [n], ({\mathbb P}^1)^{\wedge^{l(j_1) + \cdots +l(j_k)}})$.
\begin{theorem} (Associativity of $\mu$)
\label{assoc}
Let $\gamma _k :NZ(E(\Sigma _k)) \otimes NZ(E(\Sigma _{n_1})) \otimes \cdots \otimes
NZ(E (\Sigma _{n_k})) \ra NZ(E (\Sigma _{\Sigma _i n_i}))$ denote the pairing 
defined by the operad-structure on $\{NZ(E(\Sigma _k))|k\}$. Then the following diagram commutes strictly:
\begin{multline}
\label{assoc.mu}
{\diagram
{NZ(E(\Sigma _k)) \otimes NZ(E(\Sigma _{j_1})) \otimes ... \otimes NZ(E(\Sigma_{j_k})) \otimes {\underline {\mathbb Z }}^{\overset j \otimes}} \stackrel{\gamma \otimes id}{\longrightarrow} \ddto_{regroup} &{NZ(E(\Sigma_j)) \otimes {\underline {\mathbb Z }}^{\overset j \otimes}} \dto^{\mu}\\
& {\underline {\mathbb Z}}\\
{NZ(E(\Sigma_k)) \otimes NZ(E(\Sigma_{j_1})) \otimes {\underline {\mathbb Z}}^{\overset {j_1} \otimes}\otimes ... \otimes NZ(E(\Sigma_{j_k})) \otimes {\underline {\mathbb Z }}^{\overset {j_k} \otimes}} \overset { id \otimes \mu ^k}{\longrightarrow}  &{NZ(E(\Sigma _k)) \otimes {\underline {\mathbb Z }}^{\overset k \otimes}} \uto_{\mu}
\enddiagram}
\end{multline} 
\end{theorem}
\begin{proof} The key observation here is that the pairing $\mu'$ as in 
~\eqref{key.pairing.4}, and ~\eqref{assoc.key.2} defined above commutes with shuffle maps. The commutativity of the above diagram follows essentially 
from this observation. However, we provide some details below, mainly for the sake of completeness.
\vskip .3cm
Let $s_{j_1}, \cdots s_{j_k}$ denote simplices in $NZ(E(\Sigma _{j_1})), \cdots, NZ(E(\Sigma _{j_k}))$ of dimensions $p_1, \cdots p_k$ respectively. Let
$t_k$ denote a $q$-simplex of $N(ZE(\Sigma _k))$. 
\vskip .3cm
We will first consider the $q+p_1 + \cdots +p_k$-simplex of $NZ(E(\Sigma _{j_1+ \cdots j_k}))$ defined by $\gamma (t_k, s_{j_1}, \cdots, s_{j_k})$: see Definition ~\ref{Barratt.Eccles.simp}. 
One first produces a $(q, p_1, \cdots, p_k)$ multi-simplex by pairing these as in Definition ~\ref{Barratt.Eccles.simp}. Next one applies a shuffle map to produce a $q+\Sigma _i p_i$ simplex from these.
The associativity of the shuffle maps entering into its definition shows that this may be obtained
as follows. One first applies a $(p_1, \cdots, p_k)$-{\rm {shuffle}} $(\mu_{j_1}, \cdots \mu_{j_k})$ to the simplex $(s_{j_1}, \cdots , s_{j_k})$
to obtain the $p_1+ \cdots +p_k$-simplex $(\mu_{j_1}(s_{j_1}), \cdots, 
\mu _{j_k}(s_{j_k}))$ in the image of $Z(E(\Sigma _{j_1})) \otimes \cdots \otimes Z(E(\Sigma _{j_k}))
\ra Z(E(\Sigma _{\Sigma _i j_i}))$. Next one applies a $(q, \Sigma _i p_i)$ -{\rm {shuffle}} $(\mu_k, \mu_j)$ to
$(t_k, (\mu_{j_1}(s_{j_1}), \cdots, \mu _{j_k}(s_{j_k})))$ to obtain 
a $q+\Sigma _i p_i$-simplex in the image of $Z(E(\Sigma _k))\otimes Z(E(\Sigma _{j_i})) \ra Z(E(\Sigma _{\Sigma _i j_i}))$. 
\vskip .3cm 
Next one starts with cycles $Z(j_i)_1, \cdots Z(j_i)_{j_i}$ in $\Gamma (U, {\underline {\mathbb Z}}(l(j_i)_1)), \cdots, \Gamma (U, {\underline {\mathbb Z}}(l(j_i)_{j_i}))$, $i=1, \cdots, k$. (i.e. $Z(j_i)_m \eps \Cf(U \times \Delta [n], ({\mathbb P}^1)^{\wedge ^{l(j_i)_m}})$.)
The cycle 
\[\mu(\gamma (t_k, (s_{j_1}, \cdots, s_{j_k})), Z(j_1)_1, \cdots, Z(j_1)_{j_1}, \cdots, Z(j_k)_1, \cdots, Z(j_k)_{j_k})\]
is obtained by applying the pairing of (~\ref{pairing.4}). Since the shuffle maps commute
with the pairing $\mu'$ (see ~\ref{pairing.1.shuffle}, (~\ref{shuffle.real}) and  ~\eqref{assoc.key.2}), it follows that the last pairing may be obtained as
\be \begin{equation}
\label{toprow}
{\rm {shuffle}} \circ \mu'(t_k, s_{j_1}, \cdots , s_{j_k}, Z(j_1)_1, \cdots, Z(j_1)_{j_1}, \cdots, Z(j_k)_1, \cdots Z(j_k)_{j_k})
\end{equation} \ee
Here $shuffle$ denotes the composition of
all the shuffle maps for passage from $\Delta[n] \times \Delta[p_1] \times \cdots \Delta [p_k] \times \Delta [q] \ra \Delta [n+\Sigma _i p_i+q]$.
\vskip .4cm
Next consider the pairing defined by the composition of maps in the left-most column, the bottom row and the bottom part of the right column applied to the same cycles as above. One may see that the resulting cycle may be obtained from
the cycles $Z(j_i) =\mu'(s_{j_i}, Z(j_i)_1, \cdots, Z(j_i)_{j_i}) \eps \Cf(U \times \Delta[n] \times \Delta[p_i], ({\mathbb P}^1)^{\wedge ^{l_{j_i}}})$   as follows. (Here $l_{j_i} = l(j_i)_1 + \cdots +l(j_i)_{j_i}$ where $Z(j_i)_m \eps \Cf(U \times \Delta [n], ({\mathbb P}^1)^{\wedge ^{l(j_i)_m}})$. 
\vskip .3cm
One composes this with a suitable shuffle map to obtain the 
class 
\[\mu(s_{j_i}, Z(j_i)_1, \cdots , Z(j_i)_{j_i}) \eps \Cf(U \times \Delta[n+p_i], 
({\mathbb P}^1)^{\wedge ^{l_{j_i}}}).\]
\vskip .3cm \noindent
Next one applies shuffle maps
corresponding to $(\mu _{j_1}, \cdots, \mu _{j_k})$  to obtain  classes in 
\vskip .3cm \noindent
$ \Cf(U \times \Delta [n+p], ({\mathbb P}^1)^{\wedge ^{l_{j_i}}}).$ (Recall $p = \Sigma _i p_i$.) 
\vskip .3cm
These will still be denoted $\mu(s_{j_i}, Z(j_i)_1, \cdots , Z(j_i)_{j_i})$ for the sake of notational 
simplicity.  Next one applies the pairing in (~\ref{pairing.4}) to these $k+1$
simplices as $\sigma$ varies among the vertices of the simplex $t_k$ to obtain the class 
\[\mu(t_k, \mu (s_{j_1}, Z(j_1)_{1}, \cdots, Z(j_1)_{j_1}), \cdots,
\mu(s_{j_k}, Z(j_k)_1, \cdots, Z(j_k)_{j_k}))\]
\[= {\rm {shuffle}}_3 \circ \mu'( t_k, {\rm {shuffle}}_2 \circ {\rm {shuffle}}_1\circ \mu'(s_{j_1}, Z(j_1)_1, \cdots, Z(j_1)_{j_1}), \]
\[\cdots, {\rm {shuffle}}_2 \circ {\rm {shuffle}} _1 \mu'( s_{j_k}, Z(j_k)_1, \cdots , Z(j_k)_{j_k})). \]
\vskip .3cm \noindent
Here the shuffle map marked
${\rm {shuffle}}_1$ (${\rm {shuffle}} _2$, ${\rm {shuffle}}_3$) corresponds to the passage from $\Delta [n] \times \Delta [p_i]$ to $\Delta [n+p_i]$ (the passage from $\Pi _i\Delta [n+p_i]$ to $\Delta [n+p]$, the passage from
$\Delta [q]\times \Delta [n+p] $ to $\Delta [n+p+q]$, \res.)
\vskip .3cm
In view of the observation in ~\ref{pairing.1.shuffle}, one may postpone
applying the shuffle maps so that one may identify the last product with
\be \begin{equation}
\label{pairing.6}
{\rm {shuffle}}  \mu'(t_k, {\rm {shuffle}}_4\mu'(s_{j_1}, Z(j_1)_1, \cdots, Z(j_1)_{j_1}), \cdots, {\rm {shuffle}}_4\mu'( s_{j_k}, Z(j_k)_1, \cdots , Z(j_k)_{j_k}))
\end{equation} \ee
\vskip .3cm \noindent
 Here $shuffle_4$ denotes the passage from $\Pi_i \Delta [p_i]$ to $\Delta [p]$, the last 
shuffle maps are the obvious remaining shuffle maps and 
\vskip .3cm
$\mu'(t_k, {\rm {shuffle}}_4\mu'(s_{j_1}, Z(j_1)_1, \cdots, Z(j_1)_{j_1}), \cdots, {\rm {shuffle}}_4\mu'( s_{j_k}, Z(j_k)_1, \cdots , Z(j_k)_{j_k})) \eps
\Cqf(U \times \Delta [n] \times \Delta [p] \times \Delta [q], ({\mathbb P}^1)^{\wedge ^{\Sigma _i l_{j_i}}})$.
\vskip .3cm
The $p$-simplex $\Delta [p]$ denotes the $p$ simplex of $NZ(E(\Sigma _{j_1})) \otimes \cdots \otimes NZ(E(\Sigma _{j_k}))$ given by the product $\sigma _{j_1}(s_{j_1}) \otimes \cdots \otimes \sigma _{j_k}(s_{j_k})$. 
\subsubsection{}
\label{pairing.2.shuffle}  
In view of the above observations the shuffle maps all commute with $\mu'$. Therefore, in view of the identification in (~\ref{assoc.key.2}) and the observation in (~\ref{pairing.5}), 
the cycle in (~\ref{pairing.6}) may be identified with 
\be \begin{multline}
{\mu'}(t_k , s_{j_1}, \cdots, s_{j_k}, {\underset {i=1} { \overset {j_1} \bullet}} Z(j_1)_i,
  \cdots , {\underset {i=1} { \overset {j_k} \bullet}}Z(j_k)_i)\\
= {\mu'}(t_k, s_{j_1}, \cdots, s_{j_k}, Z(j_1)_1, \cdots, Z(j_1)_{j_1}, \cdots Z(j_k)_1, \cdots, 
Z(j_k)_{j_k})
\end{multline} \ee
\vskip .3cm \noindent
followed by appropriate shuffle maps. An inspection of the definition of ~\eqref{toprow} shows this clearly identifies with the product in (~\ref{toprow}) thereby proving the theorem.
\end{proof}
\begin{theorem}
\label{alg.struct}
The motivic complexes 
${\underline {\mathbb Z}} ={\underset r \oplus}{\underline {\mathbb Z}}(r)$ and 
${\underline {\mathbb Z}}/l^{\nu} ={\underset r \oplus}{\underline {\mathbb Z}}/l^{\nu}(r)$ (for any fixed prime $l$ and $\nu >0$)
are $E^{\infty}$-algebras over the $E^{\infty}$-operad $\{NZ(E(\Sigma _k))|k\}$  in the category $\oC(({\rm {Sm}}/k)_{Nis})$. Similarly the motivic complex ${\underline {\mathbb Z}}^{et}/l^{\nu} = {\underset r \oplus}{\underline {\mathbb Z}}^{et}/l^{\nu}(r)$ is an $E^{\infty}$-algebra over the $E^{\infty}$-operad $\{NZ(E(\Sigma _k))|k\}$  in the category $\oC(({\rm {Sm}}/k)_{et})$ for each fixed prime $l$ and $\nu > 0$.\end{theorem}
\begin{proof} We will only consider explicitly  the statements for the integral motivic complex,
since the
corresponding statements for the mod-$l^{\nu}$ motivic complexes follow along the same lines. Theorem ~\ref{assoc} proves that the pairing $\mu$ satisfies the condition
 ~\eqref{operad.4} in the appendix. The fact that the pairing $\mu$ also satisfies the condition
~\eqref{operad.5} should be clear. The observations in ~\ref{mu.comm} show that it also satisfies
the condition ~\eqref{operad.algebra}. \end{proof}
\vskip .3cm
\begin{definition}
(i) ${\underline {\mathbb Z}}^{mot}$ will denote the motivic complex viewed as an algebra
over the operad $\{NZ(E(\Sigma _k))|k\}$.
\vskip .3cm 
\label{W}
(ii) Recall from \cite[Theorem 1.4, Chapter II]{km}, that there exists a functor $W$ that
converts any $E^{\infty}$-algebra tensored with $\Q$ to a quasi-isomorphic strictly commutative
differential graded algebra. We will let ${\underline \Q}^{mot}$ denote 
$W({\underline {\mathbb Z}}^{mot} {\underset {\mathbb Z} \otimes} \Q)$. 
We also let ${\underline \Q}^{mot}(n) =$ the corresponding piece of weight $n$.  
If $X$ is a smooth separated scheme of finite type over $k$, $\Q^{mot}_X = R\Gamma (X, {\underline \Q}_X^{mot})$ will be called the {\it motivic DGA} associated to $X$.\end{definition}
Next we obtain the following corollary.
\begin{corollary} ${\underline \Q}^{mot}_X$ is quasi-isomorphic
to ${\underline {\mathbb Z}}_X{\underset {\mathbb Z} \otimes} \Q$ as
a sheaf of differential graded algebras.  Moreover,
${\underline \Q}^{mot}_X(n)$ is quasi-isomorphic
to ${\underline {\mathbb Z}}_X(n){\underset {\mathbb Z} \otimes} \Q$ and if $X$
is a smooth quasi-projective variety, $\Q_X^{mot}(n)$ is quasi-isomorphic to 
$\Gamma (X, {\underline {\mathbb Z}}_X^{mot}(n)){\underset {\mathbb Z} \otimes} \Q$ for each $n \ge 0$.
\end{corollary}
\begin{proof} It follows from the proof \cite[Corollary (1.5), Part II]{km} that there is
a natural map ${\underline {\mathbb Z}}_X(n){\underset {\mathbb Z} \otimes} \Q
\ra {\underline \Q}_X^{mot}(n)$ which is a quasi-isomorphism compatible with the pairings. The quasi-isomorphism in the last statement now follows from the
observation that the integral motivic complex, being quasi-isomorphic to 
the higher cycle complex, satisfies the localization property and hence
has cohomological descent on the Zariski site of smooth quasi-projective
schemes.
\end{proof}
\begin{corollary} For any  scheme $X$ over $k$, 
${\mathbb Z}_X^{mot} = \Gamma (X, {\underline {\mathbb Z}}^{mot})$, ${\mathbb Z}^{mot}/l_X^{\nu}= \Gamma (X, {\underline {\mathbb Z}}^{mot}/l^{\nu})$ 
and ${\mathbb Z}^{et}/l_X^{\nu} =  \Gamma (X, {\underline {\mathbb Z}}^{et}/{l^{\nu}})$ are  algebras over the operad
$NZ(E\Sigma)$. Moreover if $X$ is quasi-projective,
\[{\mathbb Z}_X^{mot} \simeq R\Gamma (X, {\underline {\mathbb Z}}^{mot}), \ {\mathbb Z}^{mot}/l_X^{\nu} \simeq  R\Gamma (X, {\underline {\mathbb Z}}^{mot}/l^{\nu}) \, and \]
\[{\mathbb Z}^{et}/l_X^{\nu} \simeq  R\Gamma (X, {\underline {\mathbb Z}}^{et}/{l^{\nu}}).\]
\vskip .3cm
For a general smooth scheme $X$ of finite type over $k$,
${\mathbb Z}_X^{mot} = R\Gamma (X, {\underline {\mathbb Z}}^{mot})$, ${\mathbb Z}^{mot}/l_X^{\nu}= R\Gamma (X, {\underline {\mathbb Z}}^{mot}/l^{\nu})$ 
and ${\mathbb Z}^{et}/l_X^{\nu} = R \Gamma (X, {\underline {\mathbb Z}}^{et}/{l^{\nu}})$ (where the derived functors are taken using a Godement
resolution on the Nisnevich site of $X$) are  algebras over
$NZ(E\Sigma) \otimes \End_{{\mathcal Z}}$ where the tensor product of operads may be defined as in \cite[4.1]{Moerd}.
\end{corollary}
\begin{proof} This follows immediately from Proposition ~\ref{norm.operad} in view of the following observations. Observe that the Nisnevich site 
of a scheme of finite type over $k$ has finite cohomological dimension. Therefore, the derived
functor $R\Gamma (X, {\underline {\mathbb Z}})$ may be computed as $Tot N(\{\Gamma (X, G^n ({\underline {\mathbb Z}}))|n \})$. This is the total complex of the double complex obtained by normalizing in the 
cosimplicial direction applied to the cosimplicial object 
$\{\Gamma (X, G^n ({\underline {\mathbb Z}}))|n \}$. Observe also that ${\underline {\mathbb Z}}^{et}/
{l^{\nu}}(r) \simeq \mu_l^{\nu}(r)[0]$ so that 
once again one may compute $R \Gamma (X_{et}, {\underline {\mathbb Z}}^{et}/{l^{\nu}}r))$ by using the
total complex construction as above. Since the Nisnevich site of a scheme of finite type over $k$ has finite cohomological dimension,
the total complex construction here is well-behaved.

\end{proof}
\begin{remarks}1. One may also observe that the above $E_{\infty}$-structures on the motivic complex ${\mathbb Z}^{mot}/l^{\nu}_X$ and 
${\mathbb Z}^{et}/l_X^{\nu}$ are compatible. To see this, it suffices to observe that the natural map ${\mathbb Z}^{mot}/l^{\nu} \ra
R\epsilon_*(\epsilon ^{-1}{\mathbb Z}^{mot}/l^{\nu}) \ra R \epsilon_*({\mathbb Z}^{et}/l^{\nu})$ is a map of $E_{\infty}$-algebras over
 the simplicial Barratt-Eccles operad $\{ NZ(E\Sigma _n))|n \ge 0\}$.
\vskip .3cm
2. It will be convenient and often necessary to obtain actions
by other operads, (for example, a more geometric form of the Barratt-Eccles operad) on the
motivic complex. These actions will be induced by the above action of
the simplicial Barratt-Eccles operad and will be produced by defining a
map from the new operad to the  simplicial Barratt-Eccles operad. It is
convenient to invoke model structures on the category of operads (provided, for example, by Hinich : see \cite{hin}) to define such maps. These will be explored elsewhere.
\end{remarks}
\vskip .3cm
\subsection{More intuition on the graph construction}
This section is included only to provide more insight into the graph construction discussed above and is not used elsewhere
in the body of the paper. {\it To motivate and clarify our construction, we will again consider first explicitly how one proves (first order) 
homotopy commutativity of the product on the motivic complexes.} 
\vskip .3cm
We will consider a construction using  the motivic complex of weight $n$ defined as  
\[{\underline {\mathbb Z}}(n) = C^*\Cqf( \quad, {\mathbb A}^{n})[-2n].\]
 We will show that while one may construct
the higher order homotopies explicitly this way, the structure that one obtains on the graded motivic complex $\oplus _{n \ge 0} {\mathbb Z}(n)$ will 
be only one of an {\it algebra over  a colored  operad.} Nevertheless this construction clarifies 
our basic ideas and therefore, we will begin by discussing it briefly. It will become clear this is a variant of the
 construction defined eralier using $({\mathbb P}^1)^{\wedge ^n}$ in the place of ${\mathbb A}^n$ 
which does indeed provide the structure of an algebra on $\oplus _{n \ge 0} {\mathbb Z}(n)$ over the Barratt-Eccles operad.
\vskip .3cm
Since ${\underline {\mathbb Z}}(n) = C^*\Cqf( \quad, {\mathbb A}^{n})[-2n]$, the product is defined by the pairing
\be \begin{equation}
\label{pairing.1}
\Cqf(\quad \times \Delta [.], {\mathbb A}^l) \otimes \Cqf(\quad \times \Delta [.], {\mathbb A}^m) 
\ra \Cqf(\quad \times \Delta [.], {\mathbb A}^{l+m})
\end{equation} \ee
of simplicial abelian sheaves, which may be described as follows.
Here $\Delta[.] $ denotes the cosimplicial scheme $\{\Delta[n]|n\}$ with the
obvious structure maps. We start with correspondences $\Gamma _1$ on $U \times \Delta[n] \times {\mathbb A}^l$ and $\Gamma_2$ on $U \times \Delta[n ] \times {\mathbb A}^m$ and take their external product to define
a correspondence $\Gamma _1 \times \Gamma _2$ on $U \times U \times \Delta [n] \times \Delta [n] \times {\mathbb A}^{l+m}$. Next we pull-back by the diagonal $\Delta: U \times \Delta[n] \ra U \times \Delta [n] \times U \times \Delta [n]$. This
defines the pairing of simplicial abelian groups given in ~\eqref{pairing.1}. Finally we may apply the normalization functor to pass from a simplicial abelian group to the associated chain complex. 
\vskip .3cm
{\it Key observation:} Alternatively one may first take the chain complexes associated to the  simplicial abelian groups 
on the left-side of ~\eqref{pairing.1} and define a pairing on them. This will involve passing to the corresponding
simplicial abelian groups using shuffle maps, making use of the  pairing in ~\eqref{pairing.1} at the level of simplicial abelian groups and
passing to the chain complex associated to the simplicial abelian group on the right-side of ~\eqref{pairing.1}. Since 
the shuffle maps from the product of the resulting chain complexes to the 
chain complex associated to the product of the simplicial abelian groups strictly commute with the action of the
symmetric group $\Sigma _2$ permuting the two factors, the need for homotopy commutativity arises only from
the switching of the two weight-factors ${\mathbb A}^l$ and ${\mathbb A}^m$. (In particular, the required homotopy may be  constructed from an ${\mathbb A}^1$-homotopy in ${\mathbb A}^{l+m}$: this observation 
should shed some perspective on the detailed constructions in ~\ref{max.ideal.1} through ~\ref{graph.const.2}.)
\vskip .3cm
\begin{remark} One may contrast the pairing ~\eqref{pairing.1} with the pairing on singular cohomology. There one starts with a pairing of cosimplicial abelian 
groups and on passage to the associated co-chain complexes, the resulting pairing involves the Alexander-Whitney maps which do {\it not } commute 
strictly with the action of the symmetric group. The failure of this strict commutativity leads to the existence of cohomology operations in singular cohomology. Thus, already one can see there are {\it important differences} between the pairing on the motivic complexes and the one on the singular co-chain complex of a topological space. \end{remark}
\subsection{}
\label{max.ideal.1} 
Assume that we are given $2$ cycles $Z_1 \eps \Cqf(U \times \Delta[n], {\mathbb A}^l)$ and
$Z_2 \eps \Cqf(U \times \Delta[n], {\mathbb A}^m)$. After pull-back by the diagonal
$\Delta: U \times \Delta[n] \ra U\times \Delta[n] \times U \times \Delta [n]$,  $Z_1 \times Z_2$
($Z_2 \times Z_1)$ defines the cycle $\Delta^*(Z_1 \times Z_2) \eps \Cqf(U \times \Delta[n], {\mathbb A}^l \times {\mathbb A}^m)$ (the cycle $\Delta^*(Z_2 \times Z_1) \eps \Cqf(U \times \Delta [n], {\mathbb A}^m \times {\mathbb A}^l)$, \res). {\it At this point we may assume that the cycles $Z_i$ are in fact closed sub-schemes of
$U \times \Delta[n] \times {\mathbb A}^l$ and $U \times \Delta [n] \times {\mathbb A}^m$}. By replacing $\Delta^*(Z_1 \times Z_2)$ with an irreducible component, if necessary, we may assume $\Delta^*(Z_1 \times Z_2)$ is a closed sub-scheme of $U \times \Delta [n] \times {\mathbb A}^{l+m}$. 
\vskip .3cm
{\it The basic idea is to attempt to define a subscheme of $U \times \Delta[n] \times {\mathbb A}^1 \times {\mathbb A}^{l+m}$ that joins a point on the scheme $\Delta^*(Z_1 \times Z_2)$ with
the corresponding point on $\Delta^*(Z_2 \times Z_1)$ by a line and show that this defines a cycle in 
$\Cqf(U \times \Delta [n] \times \Delta [1], {\mathbb A}^{l+m})$. This will provide the first order homotopy relating the classes 
$\Delta^*(Z_1 \times Z_2)$ and $\Delta^*(Z_2 \times Z_1)$.} Since there are technical difficulties for constructing
this as a correspondence in $U \times \Delta[n] \times {\mathbb A}^1 \times {\mathbb A}^{l+m}$ as explained below in 
Examples ~\ref{countereg}, we carry out this construction
in the scheme, $U\times \Delta [n] \times \Gamma _{\phi_{s_1}}$ (as defined below) which is naturally isomorphic to the scheme 
$U \times \Delta[n] \times {\mathbb A}^1 \times {\mathbb A}^{l+m}$.
\vskip .3cm 
 \subsection{The graph construction}
\label{graph.constr.1}
 The detailed construction is as follows. There is an obvious action by the group $\Sigma_2$ on $\sqcup_{l, m}{\mathbb A}^{l+m} = \sqcup_{l, m}{\mathbb A}^l \times {\mathbb A}^m$ switching the two factors ${\mathbb A}^l$ and ${\mathbb A}^m$. If $\sigma \eps \Sigma _2$
is the non-identity element, $\sigma (x_1, \cdots x_l, y_1, \cdots , y_m) = (y_1, \cdots, y_m, x_1, \cdots, x_l)$. Let $s_1 =(id, \sigma)$ denote the
obvious $1$-simplex of $E\Sigma _2$. We define a map of schemes $\phi_{s_1}:  {\mathbb A}^1 \times  (\sqcup _{l, m \ge 1} {\mathbb A}^{l}  \times {\mathbb A}^{m}) \ra  \sqcup_{l, m\ge 1}({\mathbb A}^{m+l})$ by 
\vskip .3cm
\[\phi_{s_1}(t, x_1, \cdots, x_l, y_1, \cdots y_m) = (1-t).(x_1, \cdots, x_l, y_1, \cdots, y_m) + t(\sigma(x_1, \cdots, x_l, y_1, \cdots,y_m)).\] 
\vskip .3cm \noindent
The restriction of $\phi_{s_1}$ to ${\mathbb A}^1 \times {\mathbb A}^l \times {\mathbb A}^m$ will be denoted $\phi_{s_1, l, m}$: if $l, m$ are fixed
throughout our discussion, as below, we will abbreviate this to $\phi_{s_1}$.
Now the {\it graph}
of $\phi_{s_1,l,m}$, $\Gamma _{\phi_{s_1,l,m}}$, is the closed subscheme of $ {\mathbb A}^1 \times ( {\mathbb A}^{l+m} \times {\mathbb A}^{l+m})$ 
defined by the ideal    
\vskip .3cm \noindent
$\{(t, x_1, \cdots, x_l, y_1, \cdots, y_m, \bar x_1, \cdots \bar x_l, \bar x_{l+1}, \cdots, \bar x_{l+m})| $
\vskip .3cm 
$ \quad \quad \quad  \bar x_i - p_i((1-t)(x_1, \cdots x_l, y_1, \cdots y_m)+ t. \sigma(x_1, \cdots, x_l, y_1, \cdots, y_m))|i \}.$
\vskip .3cm \noindent
(Here $p_i$ denotes projection to the $i$-th factor.)
\subsubsection{}
\label{graph.const.2}
Moreover, the graph of $\phi_{s_1,l,m}$ is isomorphic (as a scheme) by projection to the domain
of $\phi_{s_1, l, m}$, i.e. to ${\mathbb A}^1 \times {\mathbb A}^{l+m}$. In view of this, it may be important to observe that $U \times \Delta [n] \times \Gamma_{\phi_{s_1, l, m}}$ is isomorphic to the product of $U \times \Delta [n] $ and  ${\mathbb A}^1 \times {\mathbb A}^{l+m}$ in the category of schemes though it is imbedded in $U \times \Delta [n] \times {\mathbb A}^1 \times {\mathbb A}^{l+m} \times {\mathbb A}^{l+m}$ not as
a co-ordinate hyperplane, i.e. not by putting some of the co-ordinates in ${\mathbb A}^{l+m} \times {\mathbb A}^{l+m}$ zero. 
\vskip .3cm
We consider a  closed integral subscheme $Z \subseteq U \times \Delta [n] \times {\mathbb A}^{l+m}$ so that its projection to $U \times \Delta [n]$ is 
quasi-finite and dominant.
For the most part we may assume $Z=\Delta ^*(Z_1 \times Z_2)$ chosen as in ~\ref{max.ideal.1}, but the graph-construction does not require this.
 Its inverse image under the obvious projection $p: U \times \Delta [n] \times {\mathbb A}^1 \times
{\mathbb A}^{l+m} \ra U \times \Delta [n] \times {\mathbb A}^{l+m}$ will define the closed subscheme
$p^{-1}(Z)$ (which $=p^{-1}(\Delta^*(Z_1 \times Z_2))$, if $Z= \Delta^*(Z_1 \times Z_2)$).
 Now the graph of the restriction of $id \times \phi_{s_1}$ to this subscheme
defines the closed subscheme $\Gamma _{id \times \phi_{s_1}|p^{-1}(Z))}$ of $\Gamma _{id_{U \times \Delta [n ]}} \times
\Gamma _{\phi_{s_1}} $ which is contained in $ U \times \Delta [n] \times U \times \Delta [n] \times {\mathbb A}^1 \times {\mathbb A}^{l+m} \times
{\mathbb A}^{l+m}$. Since the composition $U \times \Delta [n] {\overset {\Delta} \ra} U \times \Delta [n] \times U \times \Delta [n] {\overset {pr_1} \ra} U \times \Delta [n]$ is the identity, and hence in particular proper, the image of 
$\Gamma _{id \times \phi_{s_1}|p^{-1}(Z))}$ under the projection $pr_1 \times id$ to $U \times \Delta [n] \times {\mathbb A}^1 \times {\mathbb A}^{l+m} \times {\mathbb A}^{l+m}$ is closed. Since it is the image of an irreducible scheme, it is also irreducible. {\it We will denote this by $s_1(Z)$.} 
We summarize some of the main properties of this construction here:
\subsubsection{}
\label{graph.const.props}
\begin{itemize}
\item
Let $\sigma(Z)$ denote the image of $Z$ under the permutation of the two factors ${\mathbb A}^l$ and 
${\mathbb A}^m$ in $U \times \Delta [n] \times {\mathbb A}^{l+m}$. (When $Z = \Delta^*(Z_1 \times Z_2)$, clearly $\sigma(Z) = \Delta^*(Z_2 \times Z_1)$.)
It is clear from the above definition that the scheme $ s_1(Z)$ consists of triples, consisting of a point of $Z$, a point $t \eps {\mathbb A}^1$ together with
a point on the  line joining this point of $Z$  with the corresponding point of
$\sigma(Z)$ parametrized by $t \eps {\mathbb A}^1$. Therefore, the projection of
$ s_1(Z)$ to $U \times \Delta [n] \times {\mathbb A}^1$ is in fact {\it quasi-finite}.
It is also {\it dominant} since the projection of $Z \ra U \times \Delta[n]$ is dominant. 
\item
$p^{-1}(Z)$ being a product of ${\mathbb A}^1$ and $Z$ is evidently  irreducible. Since the graph $\Gamma _{id \times \phi_{s_1}|p^{-1}(Z)}$ is isomorphic to $p^{-1}(Z))$, it is also {\it irreducible}. Therefore, so is $ s_1(Z))$.
\item 
The isomorphism of the graph $\Gamma _{\phi_{s_1,l,m}}$ with 
${\mathbb A}^1 \times {\mathbb A}^{l+m}$ shows the following: let
$\Cqf(U \times \Delta[n] \times \Gamma _{\phi_{s_1}})= $ all correspondences on $U \times \Delta[n] \times \Gamma _{\phi_{s_1}}$ whose projection to $U \times \Delta [n] \times {\mathbb A}^1$ (where ${\mathbb A}^1$ is contained in the  domain of $\phi_{s_1}$) is quasi-finite and dominant. Then this abelian group identifies with 
$ \Cqf(U \times \Delta [n] \times {\mathbb A}^1, {\mathbb A}^{l+m})$ = the correspondences on $U \times \Delta [n] \times {\mathbb A}^1 \times {\mathbb A}^{l+m}$ whose projections to $U \times \Delta [n] \times {\mathbb A}^1$ are quasi-finite and dominant. {\it Since this isomorphism is natural in $U \times \Delta [n]$, by varying $U \times \Delta [n]$, we obtain the identification 
\be \begin{equation}
\label{key.ident}
\Cqf(\  \times \Delta [.] \times {\mathbb A}^1, {\mathbb A}^{l+m}) \cong \Cqf(\  \times \Delta[.] \times \Gamma _{\phi_{s_1}})
\end{equation} \ee
\vskip .3cm \noindent 
of  the corresponding simplicial abelian sheaves}. In view of this isomorphism, we will identify the associated (co-)chain complex of the term on the right with the associated  (co-)chain complex of the term on the left, i.e.
$\Gamma (U \times {\mathbb A}^1, {\mathbb Z}(l+m)$. {\it Henceforth we will denote the term ${\mathbb A}^1$ appearing above by $\Delta [1]$ to signify that it is {\it not} a weight-factor.} 
\item 
Next assume $Z= \Delta^*(Z_1 \times Z_2)$ as in ~\ref{max.ideal.1}. Then  $s_1(\Delta ^*(Z_1 \times Z_2))$ defines a class
in $\Cqf(U \times \Delta [n] \times  \Gamma _{\phi_{s_1}}) \cong \Cqf(U \times \Delta [n] \times {\mathbb A}^1, {\mathbb A}^{l+m})$.  Restricting to the two faces of $\Delta[1] \subseteq \Gamma_{\phi_{s_1}}$ provides the two classes 
in $\Cqf(U \times \Delta [n] \times \Gamma _{id})$ and $\Cqf(U \times \Delta [n] \times \Gamma _{\sigma})$: {\it the restriction maps 
\[\Cqf(\  \times \Delta[.] \times \Gamma _{\phi_{s_1}}) \ra  \Cqf(\  \times \Delta[.] \times \Gamma _{\tau})\]
\vskip .3cm \noindent
are both quasi-isomorphisms in view of the homotopy property of the motivic complexes and the isomorphism in ~\eqref{key.ident}.} (For each permutation
$\tau \eps \Sigma _2$, we let $\Gamma _{\tau}$ denote the graph of the induced map $\tau^*:{\mathbb A}^{l+m} \ra {\mathbb A}^{l+m}$.
$\Cqf(U \times \Delta [n] \times \Gamma _{\tau})$ denotes the correspondences on $U \times \Delta[n] \times \Gamma _{\tau}$ whose projections to $U \times \Delta [n] $ are quasi-finite and dominant.)
One may further take the image of these classes under the composite map $\Gamma _{\tau} \subseteq {\mathbb A}^{l+m} \times {\mathbb A}^{l+m} {\overset {pr_2} \ra} {\mathbb A}^{l+m}$ (where the last map is the projection to the second factor) to obtain classes that identify with
 $\Delta^*(Z_1 \times Z_2))$ and
$\sigma ^*(\Delta ^*(Z_1 \times Z_2))$, \res. (One may also see readily that the above classes in $\Cqf(U \times \Delta [n] \times \Gamma _{\tau})$ map isomorphically to their images in $\Cqf(U \times \Delta [n], {\mathbb A}^{l+m})$ under the  projection ${\mathbb A}^{l+m} \times {\mathbb A}^{l+m} {\overset {pr_2} \ra} {\mathbb A}^{l+m}$.)
\end{itemize}
\vskip .3cm \noindent
Observe that, in the description of the first order homotopy as above, we only considered the
$1$-simplices of $E\Sigma_2$ of the form $s_1=(id, \sigma)$, $\sigma \eps \Sigma_2$. 
We may extend the above construction to define an action by all $1$-simplices of $E\Sigma_2$ which are of the form
$s_1=(\sigma_0, \sigma_1)$, $\sigma_i \eps \Sigma _2$ as follows: we simply replace the subscheme
$\Delta^*(Z_1 \times Z_2)$ by $\sigma_0^*(\Delta^*(Z_1 \times Z_2))$ where $\sigma_0^*$ is the 
map induced by the permutation action $\sigma_0: \sqcup_{l,m}{\mathbb A}^l \times {\mathbb A}^m \ra \sqcup_{l,m} {\mathbb A}^l \times {\mathbb A}^m$ and apply the same constructions as before
with $\sigma_1$ playing the role of $\sigma$. 
Therefore, the above construction provides a pairing:
\vskip.3cm
$Z(\{1-simplices \, s_1 =(\sigma_0, \sigma_1) \eps E\Sigma_2\}) \otimes \Cqf(U \times \Delta [n], {\mathbb A}^l) \otimes \Cqf(U \times \Delta [n], {\mathbb A}^m) \ra \Cqf(U \times \Delta [n] \times \Gamma_{\phi_{s_1}}) \cong \Cqf(U \times {\mathbb A}^1 \times \Delta[n] , {\mathbb A}^{l+m}) $
\vskip .3cm \noindent
where $Z(\{1-simplices \, s_1 =(\sigma_0, \sigma_1)  \eps E\Sigma_2 \})$ is the free abelian group on the $1$-simplices of $E\Sigma _2$.
This pairing is clearly compatible with restriction to the faces of the $n$-simplex $\Delta [n]$. 
By defining the graph $\Gamma_{\phi_{\sigma_0}}$ associated to a $0$-simplex $\sigma _0 \eps E\Sigma _2$ as the graph of the permutation 
$\sigma_0$ applied to $\sqcup _{l, m}{\mathbb A}^l \times {\mathbb A}^m$, one extends the last pairing to a pairing
\be \begin{equation}
\label{key.pairing.0}
{\underset s \oplus} Z(s) \otimes \{\Cqf(U \times \Delta[n], {\mathbb A}^l)|n\} \otimes \{\Cqf(U \times \Delta [n], {\mathbb A}^m)|n\} \ra \{\Cqf(U \times \Delta [n] \times \Gamma _{\phi_s})|n\}
\end{equation} \ee
\vskip .3cm \noindent
where the direct sum is taken over all the $1$-simplices in $E\Sigma _2$ and $Z(s)$ 
denotes the chain complex obtained by first taking the (free) simplicial abelian group 
on the $1$-simplex $s$ and then by normalizing it.
 Moreover, restricting to the two faces of $\Delta[1] \subseteq \Gamma_{\phi_{s_1}}$ provides the two classes corresponding to 
$\sigma_i^*(\Delta^*(Z_1 \times Z_2))$, $i=0,1 $ as observed above.  
\vskip .3cm 
Observe that both the left and right-sides of the  pairing in ~\eqref{key.pairing.0} are double-complexes whose bi-degrees are indexed by $n$ and the degree of terms in the complex $Z(s)$. The pairing is compatible with restrictions 
to the faces of $\Delta [n]$ as well as to the faces of the $1$-simplices $ s_1 \eps E\Sigma _2$ as the above descriptions show. The descriptions above show that this is compatible with the differentials of the complexes on either side. {\it Therefore, this is indeed a pairing of double complexes.}  Moreover, if $\tau \eps \Sigma_2$, it acts on a simplex $s_1=(\sigma _0, \sigma _1)$ by $\tau \circ s_1 = (\tau\circ \sigma _0, \tau \circ \sigma _1)$. The resulting pairing lands
in $\Cqf(U \times \Delta [.] \times \Gamma _{\phi_{\tau \circ s_1}})$: one may see that this is the same as the cycle obtained by applying $\tau$ first to $(Z_1, Z_2)$ and then applying the pairing with $\Gamma_{\phi_{s_1}}$. This completes the construction of an {\it explicit first order homotopy for the pairing of the motivic complexes considered in ~\eqref{pairing.1}}.
 \vskip .3cm
 It is important to observe that the target of this pairing is the complex $\{\Cqf(\quad  \times \Delta [.] \times \Gamma _{\phi_s})|n\}$, {\it i.e. it varies depending on the $1$-simplex $s \eps E\Sigma _2$}. 
It is possible to make an identification 
\[\{\Cqf(\quad \times \Delta [.] \times \Delta [1], {\mathbb A}^{l+m})|n\} {\overset {\cong} \ra } \{\Cqf(\quad  \times \Delta [.] \times \Gamma _{\phi_s})|n\}\]
 of complexes of abelian sheaves, so that the 
pairing ~\eqref{key.pairing.0} may be viewed as an {\it avatar} of a pairing:
\[NZ(sk_1E(\Sigma _2)) \otimes {\underline {\mathbb Z}}(l)  \otimes  {\underline {\mathbb Z}}(m) \ra {\underline {\mathbb Z}}(l+ m).\]
\vskip .3cm \noindent
(Recall that the complex $NZ(sk_1E(\Sigma _2))$ is trivial in degrees greater than 1.)
However, the above identification clearly loses the extra information in the complex $ \{\Cqf(\quad  \times \Delta [.] \times \Gamma _{\phi_s})|n\}$, coming in particular from the graph $\Gamma _{\phi_s}$. Keeping this complex as the target of the
above pairing (and on considering the iterated pairing involving higher dimensional simplices of $E\Sigma _k$), it becomes
necessary to replace the operad $\{Z(E \Sigma _k)|k\}$ with a related {\it colored operad} (where the $n$-simplices (and the $(n_1, \cdots, n_k)$-multi-simplices) are the colors) and {\it the last pairing
as one between  such a colored operad and a chain complex.} (Colored operads and algebras over such operads have just begun to appear in the literature: see \cite{BM2} and \cite{Lein}.)
  Rather than adopt this approach, we modify the graph construction below using projective spaces in the place
of the affine space ${\mathbb A}^n$ and obtain an action of the Barratt-Eccles operad itself on the motivic complexes. The above discussion is put in here mainly to motivate the constructions below and to point out the intricacies of our
constructions. 
\vskip .3cm 
Since the isomorphism 
 $\Cqf(\  \times \Delta [.] \times {\mathbb A}^1, {\mathbb A}^{l+m}) \cong \Cqf(\  \times \Delta[.] \times \Gamma _{\phi_{s_1}})$ 
is compatible with restriction to the faces of $\Delta [n]$ and also to $\{0, 1\} \subseteq {\mathbb A}^1 = \Delta [1]$, this construction
{\it does provide an explicit first order homotopy for the pairing  of motivic complexes considered in} ~\eqref{pairing.1}.
\begin{examples} 
\label{countereg} 
In the graph construction ~\ref{graph.constr.1} and the ensuing discussion, one may be tempted to replace the graph $\Gamma _{\phi_{s_1}}$ by the scheme ${\mathbb A}^1 \times {\mathbb A}^{l+m}$, where the ${\mathbb A}^{l+m}$ is the target of the map $\phi_{s_1}$.
This would mean one will need to replace the scheme $s_1(\Delta ^*(Z_1 \times Z_2))$ constructed in ~\ref{graph.const.2} by its image under the composite map
$ s_1(\Delta ^*(Z_1 \times Z_2)) \subseteq  U \times \Delta [n] \times \Gamma _{\phi_{s_1}} \subseteq U \times \Delta [n] \times \Delta [1] \times {\mathbb A}^{l+m} \times {\mathbb A}^{l+m} \ra  U \times \Delta [n] \times \Delta [1] \times {\mathbb A}^{l+m}$, where the last map is dropping the
first-factor of ${\mathbb A}^{l+m}$. This would be essentially projecting to the image of $\phi_{s_1}(p^{-1}(\Delta ^*(Z_1 \times Z_2))$. However,
the following counter-examples show this image (which will be denoted $Y$ in the following examples) may not be closed in $U \times \Delta [n] \times \Delta [1] \times {\mathbb A}^{l+m}$ and its closure may not be quasi-finite over $U \times \Delta [n] \times \Delta [1]$. 
\begin{enumerate}
\item[(1)]
Let $U = Spec \, k$, $n=l=m=1$. Then let $Z=\{(x, x, y) \eps {\mathbb A}^1 \times {\mathbb A}^1 \times {\mathbb A}^1| xy=1\}$. Now the projection $Z \ra U \times {\mathbb A}^1$ (where ${\mathbb A}^1$ denotes the first ${\mathbb A}^1$) is quasi-finite. 
 (In fact, this projection  is {\it not finite}, but only quasi-finite.)
Now consider 
\[Y=\{(x, t, tx+(1-t)y, (1-t)x+ty)|xy=1\}:\]
 this is the image of  $\phi_{s_1}(p^{-1}(Z))$ contained in $U \times {\mathbb A}^1 \times {\mathbb A}^1 \times {\mathbb A}^1 \times {\mathbb A}^1$. For $t=1/2$, the fiber of $Y$ over $t$ will be denoted 
$Y_{1/2}$. Now $Y_{1/2} = \{(x, 1/2, (1/2)x+(1/2)y, (1/2)y+(1/2)x)|xy=1\} =\{(x, 1/2, (1/2)x +(1/2x), (1/2x) +(1/2)x)|x \ne 0\}$. There is no limit as $x\ra 0$, so that in this example, $Y_{1/2}$ and $Y$ are closed in $U \times {\mathbb A}^1 \times {\mathbb A}^1 \times {\mathbb A}^1 \times {\mathbb A}^1$. 
\vskip .3cm
However, one can modify this example a bit to get another example, where the corresponding $Y$ will not be closed. Here we take $Z=\{(x, x, y, x, z)|xy=1, xz =-1\}$
viewed as a closed subscheme of ${\mathbb A}^1 \times {\mathbb A}^2 \times {\mathbb A}^2$, i.e. $U = Spec \, k$ and $n=1, l=m=2$.  
 Now the projection $Z \ra U \times {\mathbb A}^1$ is quasi-finite.
 Then 
\[Y=\{(x, t, t(x,y)+(1-t)(x,z), (1-t)(x,y)+t(x,z)|xy=1, xz=-1\}\] which is 
the image of  $\phi_{s_1}(p^{-1}(Z))$ contained in $U \times {\mathbb A}^1 \times {\mathbb A}^1 \times {\mathbb A}^2 \times {\mathbb A}^2$.
 Therefore, 
 $Y_{1/2}=\{(x, 1/2, ((1/2)x+(1/2)x), ((1/2)y+(1/2)z), ((1/2)x+(1/2)x), ((1/2)z+(1/2)y))|xy=1, xz=-1\}=\{(x, 1/2, x, 0, x, 0)|x \ne 0 \}$. Clearly this has $(0, 1/2, 0, 0, 0, 0)$ as a limit point which is outside of $Y$. 
Therefore, $Y_{1/2}$ and $Y$ are not closed.
In this case though, the closure of $Y$ just adds the point $(0, 1/2, 0, 0, 0, 0)$ so that this closure is still quasi-finite over the product of the first two ${\mathbb A}^1$.
\item [(2)]
Next we will construct an example, where $Y$ is not closed and the closure of $Y$ is not quasi-finite over 
$U \times \Delta [n] \times \Delta [1]$.
Here again, $U= Spec \, k$ and $n=2$ and $l=m=1$. This example is obtained from the blow-up of ${\mathbb A}^2$ at the origin with a divisor at infinity removed,
so that we obtain an affine scheme over ${\mathbb A}^2$.
Let $Z=\{(x_1, x_2, y_1, z_1)|x_1y_1=1, x_1z_1 =x_2\}$. (Recall the blow-up of ${\mathbb A}^2$ at the origin is given by equations
$x_1u_2-x_2u_1=0$ in ${\mathbb A}^2 \times {\mathbb P}^1$, where $(x_1, x_2)$ are parameters for ${\mathbb A}^2$ and $[u_1:u_2]$ are the homogeneous coordinates for ${\mathbb P}^1$. So we are letting $z_1=u_2/u_1$ and throwing out the part $u_1=0$.) Now $Z$ is clearly closed
in ${\mathbb A}^4$, as it is given by the equations $x_1y_1=1, x_1z_1 =x_2$. But the projection of $Z$ into the three factors dropping $y_1$ is not closed
in ${\mathbb A}^3$: call this $Z'$.  In fact  the projection of $Z'$ to the $(x_1, x_2)$-coordinates is not quasi-finite as the fiber over $(0, 0)$ will be a whole ${\mathbb A}^1$. 
\vskip .3cm
Now one can modify the above example to obtain a counter example where $Y$ (which is the image of  $\phi_{s_1}(p^{-1}(Z))$) is not closed and its closure is not quasi-finite for the projection to the first two factors.
Let $Z=\{(x_1, x_2, (y_1, z_1), (y_2, z_2))|x_1y_1=1, x_1y_2=-1, x_1z_1=x_2, x_1z_2=x_2\}$ viewed as a closed sub-scheme of 
${\mathbb A}^2 \times {\mathbb A}^2 \times{\mathbb A}^2$, i.e. $U = Spec \, k$, $n=2, l=m=2$. Clearly the projection of $Z$ to the $(x_1, x_2)$ 
coordinates is quasi-finite. Now 
\vskip .3cm
$Y=\{(x_1, x_2, t, t(y_1, z_1)+(1-t)(y_2, z_2), (1-t)(y_1, z_1)+t(y_2, z_2)) | $
\vskip .3cm
$ \quad \quad \quad x_1y_1=1, x_1y_2=-1, x_1z_1=x_2, x_1z_2 =x_2\}$
\vskip .3cm \noindent
 so that 
$Y_{1/2}=\{(x_1, x_2, 1/2, 1/2(y_1+y_2), 1/2(z_1+z_2), 1/2(y_1+y_2), 1/2(z_1+z_2))|x_1y_1=1, x_1y_2 =-1, x_1z_1=x_2, x_1 z_2 =x_2\}$. Clearly this
equals $\{(x_1, x_2, 1/2, 0, x_2/x_1, 0, x_2/x_1)|x_1 \ne 0\}$. Clearly this is not closed and the fiber of the closure over $x_1=0, x_2=0$ has a whole ${\mathbb A}^1$ in the
fifth and last coordinate, so that the projection of the closure of $Y_{1/2}$ to the first factor ${\mathbb A}^2$ is not quasi-finite.
\end{enumerate}
\end{examples}
The above examples make it necessary to make use of the graph $\Gamma _{\phi_{s_1}}$ and adopt the construction we have used above, i.e. if we use the definition of the motivic complex as 
${\underline {\mathbb Z}}(n) = C^*\Cqf( \quad, {\mathbb A}^{n})[-2n]$. One way to  
avoid using the graph  $\Gamma _{\phi_{s_1}}$ is to replace it by its projection to the co-domain:
however, for this to work, one needs to replace the affine spaces ${\mathbb A}^n$ appearing with $({\mathbb P}^1)^n$ and with 
${\underline {\mathbb Z}}(n) = C^*(\Cf(\quad, {\mathbb P}^{\wedge n}))[-2n]$.

\section{\bf Mixed Tate motives for smooth linear schemes over a field $k$}
The results of this section generalize the constructions of \cite{bl2}, \cite{blk} and \cite{km} for the
category of mixed Tate motives over a field. The existence of the motivic dga extends these constructions to any smooth quasi-projective scheme if one assumes the 
Beilinson-Soul\'e vanishing conjecture holds for the rational motivic cohomology of that scheme.
{\it In particular we verify this for a large class of quasi-projective smooth varieties including all projective smooth toric and spherical varieties over number fields. The main result 
is Theorem  \ref{tatemotives}. }
\vskip .3cm 
We fix a smooth quasi-projective scheme $X$
over $k$. We let $A={\Q}_X^{mot} $. We may assume therefore that $A$ is graded with
$A(r)$ denoting the part in grade $r$, where $r \ge 0$ and $A^q(r)$ denotes the part of the complex $A(r)$ in degree $q$, where $q$ is any integer. 
(Recall that $A$ has an augmentation $A \ra \Q[0]$.) Let ${\mathbb D}_{-}(A)$ denote the derived category of cohomologically bounded below $A$-modules, i.e.
differential graded $A$-modules $M = {\underset {r} \oplus} M(r)$ where $M^q(r) = (M(r))^q$ may be
non-zero for any pair of integers $(q, r)$ and so that ${\mathcal H}^q(M)(r) =0$ for all sufficiently small $q$. 
(One may first show that this derived category is equivalent to that of cohomologically bounded below {\it cell} $A$-modules in the sense of 
\cite[Part III]{km}. By construction, every cell $A$-module is {\it flat} over $A$, in the sense that the tensor product $- {\underset {A} \otimes}M $ preserves
distinguished triangles in the first argument for every cell $A$-module $M$. Then the following derived tensor product may be replaced by a tensor product.)
\vskip .3cm
Now one may define a functor 
\[Q: {\mathbb D}_{-}(A) \ra D(\mbox{${\mathbb Q}$-vector spaces}) \mbox{ by }\]
\[ Q(M) =
\Q{\overset L {\underset {A} \otimes}}M = \mbox{ the ${\mathbb Q}$-vector space of {\it indecomposable elements} of M}.\] 
See, for example, \cite[Part IV]{km}. Here $D(\mbox{$\Q$-vector spaces})$ denotes the
derived category of bounded below complexes of $\Q$-vector spaces. Observe that this category
has a natural $t$-structure, the heart of which is given by the complexes that have
cohomology trivial in all degrees except $0$. We let ${\mathcal H}_A$ denote the full
sub-category of ${\mathbb D}_A$ consisting of complexes $K$ so that ${\mathcal H}^q(Q(K)) = 0$
for all $q \ne 0$. Let ${\mathcal F}{\mathcal H}_A$ denote the full sub-category of
${\mathcal H}_A$ consisting of complexes $K$ so that ${\mathcal H}^0(Q(K))$ is a finite
dimensional $\Q$-vector space. We will make the following
assumption throughout:
\subsubsection{}
\label{hyp.dga.1}
{\it the DGA $A$ is connected in the following sense: $H^i(A)(r) = 0$ for $i<0$, $H^0(A)(r) =0$ if $ r \neq 0$ and
$H^0(A)(0) = \Q$.}
\vskip .3cm
 Now we obtain the following theorem as in \cite{km}.
\begin{theorem}
The triangulated category ${\mathbb D}_{-}(A)$ admits a $t$-structure whose
heart is ${\mathcal H}_A$. Moreover ${\mathcal F}{\mathcal H}_A$ is a graded neutral
Tannakian category over $\Q$ with fiber functor $w = {\mathcal H}^0 \circ Q$.
\end{theorem}
\begin{proof} The proof is essentially in \cite[Theorem 1.1, Part IV]{km} . (The key idea here is to use the theory of minimal models.)\end{proof}
One may apply the bar construction (see \cite[Part IV, section 1]{km}) to the algebra $A$: we will denote this
by $\bar B A$. Let $IA$ denote the augmentation ideal of $A$. We let ${\chi}_A = H^0(\bar B A)$. This is a commutative Hopf-algebra and, as in \cite[Part IV, section 1]{km} , is a polynomial algebra with its
$k$-module of indecomposable elements a co-Lie algebra which is denoted $\gamma _A$.
Now we obtain the following result.
\begin{theorem}
\label{Mixed.Tate.1}
(See \cite[Part IV, Theorem 1.2]{km}.) Assume the hypothesis (~\ref{hyp.dga.1}).
Then the following categories are equivalent:
\vskip .3cm \noindent
(i) The heart ${\mathcal H}_A$ of ${\mathbb D}_{-}(A)$
\vskip .3cm \noindent
(ii) The category of generalized nilpotent representations of the co-Lie algebra $\gamma _A$
\vskip .3cm \noindent
(iii) The category of co-modules over the Hopf-algebra $\chi_A$
\vskip .3cm \noindent
(iv) The category ${\mathcal T}_A$ of generalized nilpotent twisting matrices in $A$
\vskip .3cm \noindent
The full sub-categories of finite dimensional objects in the categories (i), (ii) and (iii)
and of finite matrices in the category (iv) are also equivalent. \end{theorem}
\begin{definition}
\label{linear}
 (Linear schemes over $k$) (i) A scheme over
$Spec \, k$ is $0$-linear if it is either empty
or isomorphic to any affine space ${\mathbb A}^n_{Spec \, k}$.
\vskip .3cm
(ii) Let $n>0$ be an integer. A scheme $Z$,  over $Spec \, k$, is $n$-linear, if there exists a
triple $(U, X, Y)$ of schemes over $Spec \, k$ so that $Y \subseteq X$ is a closed immersion
with $U$ its complement, $Y$ and one of the schemes $U$ or $X$ is $(n-1)$-linear and $Z$ is the other member in $\{U, X\}$. We say $Z$ is linear if it is
$n$-linear for some $n \ge 0$.
\vskip .3cm
(iii) Recall any reduced scheme $X$ of finite type over $Spec \, k$ is called a variety. Linear varieties over $k$ are varieties over $Spec \, k$ that are linear schemes.
\end{definition}
\begin{example}
\label{egs}
 The following are common examples of linear varieties. In these examples we fix a separably closed base field $k$ and consider only varieties over $k$.
\begin{itemize}
\item All toric varieties
\item All spherical varieties (A variety $X$ is spherical if there exists a
reductive group $G$ acting on $X$ so that there exists a Borel subgroup
having a dense orbit.)
\item
Any variety on which a connected solvable group acts with finitely many
orbits. (For example projective spaces and flag varieties.)
\item
Any variety that has a stratification into strata each of which is the
product of a torus with an affine space.
\end{itemize}
\end{example}
\vskip .3cm
\begin{remark} If the field is not separably closed, not all tori are split; therefore
the varieties appearing above need not be linear in the sense of the definition ~\ref{linear}. Over non-separably closed fields, any of the examples above will be linear if and
only if the tori appearing in the strata are all split. \end{remark}
\begin{corollary}
\label{Mixed.Tate.2}
Let $X$ denote a smooth connected projective linear variety
over a field $k$ ({\it not} necessarily separably closed), {\it or} any one of the schemes appearing in the
examples above which are also connected, projective and smooth. Assume that the Beilinson-Soul\'e conjecture holds
for the rational motivic cohomology of $Spec \, k$, i.e. $H^i_{\M}(Spec \, k; \Q (r)) = 0$ if $i<0$, $H^0_{
\M}(Spec \, k; \Q(r)) = 0$ if $r \neq 0$ and
 $H^0_{\M}(Spec \, k; \Q(0)) \cong \Q$. (For example, the above hypothesis holds if $k$ is a number field). Then the conclusions of theorem
(~\ref{Mixed.Tate.1}) hold for $X$ (i.e. with $A = \Q_X^{mot}$). Let $U = X- Y$, where $X$ and $Y$ are either projective smooth linear varieties or
any of the projective smooth schemes  appearing in the above list and that, in either case, $Y$ is closed in $X$.
Then the conclusions of Theorem ~\ref{Mixed.Tate.1} also hold for $U$. \end{corollary}
\begin{proof} It suffices to show that the DGA, $A$ appearing in the theorem is connected. If the field $k$ is not separably closed, one may find a finite separable extension $k'$
of $k$ so that all the tori in the stratification of $X$ split. By a transfer
argument, one may therefore readily reduce to the case where $k$ is separably closed. Next we will consider
the case where the scheme is {\it projective}.
In this case, the variety in question is also linear; therefore we may invoke the strong Kunneth decomposition for the class of the diagonal
$\Delta$ in $CH^*(X \times X)$. (See
\cite{josh1}.) i.e.
\be \begin{equation}
\label{s.kunneth}
\Delta = \Sigma _i \alpha _i \times \beta _i = \Sigma _i p_1^*(\alpha _i) \circ p_2^*(\beta _i)
\end{equation} \ee
where $p_i: X \times X \ra X$ is the projection to the $i$-factor, 
$\circ$ denotes the intersection product and $\alpha_i, \beta_i \eps CH^*(X)$. Now we proceed to show that any class $x \eps CH^*(X, n)$ may be written
as a linear combination
\be \begin{equation}
\label{x.1}
x = \Sigma _i \alpha _i \circ p_{1*}(p_2^*(\beta _i \circ x)) = \Sigma _i \alpha _i \circ {p_2'}^*({p_1}'_*(\beta _i \circ x))
\end{equation} 
Here $p_i': X \ra Spec \, k$ is the obvious projection.
To obtain (~\ref{x.1}), first observe that $x = p_{1*} ([\Delta ] \circ p_2^*(x))$. By the projection formula and the observation that the class $[\Delta] = \Delta _*(1)$, $1 = [X] \eps
CH^*(X ; )$, we obtain equality of the classes $[\Delta ]\circ p_2^*(x) = \Delta _* (\Delta ^*(p_2^*(x))$. Therefore
$p_{1*} ([\Delta ]\circ p_2^*(x)) = p_{1*}(\Delta_*(\Delta ^*(p_2^*(x)))) = (p_1 \circ \Delta)_* ((p_2 \circ \Delta)^*(x)) = x$. 
Now substitute the formula for $[\Delta]$ from (~\ref{s.kunneth}) and use the projection
formula to obtain the first equality in (~\ref{x.1}). The equality of this with the right-hand-side follows by flat-base-change. Observe that $\alpha _i \eps
CH^*(X, 0)$. Therefore, the hypothesis that $H^i_{\M}(Spec \, k; \Q(r)) =0$ for $i<0$ shows readily that $H^q_{\M}(X; \Q(r)) = CH^{r}(X, 2r-q; \Q ) =0$ 
for $q<0$. (In more detail: suppose $x \eps H^q_{\M}(X; \Q(r))$ for $q <0$. 
Let $\alpha _i \eps
CH^s(X, 0; \Q) = H_{\mathcal M}^{2s}(X; \Q(s))$ for some $s\ge 0$. Then if $d=dim_k(X)$, $\beta_i \eps H_{\mathcal M}^{2d-2s}(X, d-s)$, $ \beta_i \circ x \eps
H_{\mathcal M}^{q+2d-2s}(X, r+d-s)$ and 
$p'_{1*}(\beta_i \circ x) \eps H^{q-2s}_{\mathcal M}(Spec \, k; \Q(r-s))$. Since $q<0$ by assumption and $s \ge 0$, $q-2s <0$ so that $p_{1*}'(\beta_i \circ x) =0$.
The last equality is from the assumption that Beilinson-Soul\'e conjecture holds for the rational motivic cohomology of $Spec \, k$. 
Therefore $x = \Sigma _i \alpha_i \circ p_2^*(p'_{1*}(\beta_i \circ x)) =0$ as well.)
\vskip .3cm
The hypothesis that
$H^0_{\M}(Spec \, k; \Q(r)) =0$ for $r \neq 0$ implies similarly that
$H^0_{\M}(X; \Q(r)) = 0$ also for $r \neq 0$. Since $X$ is connected, the hypothesis that $H^0_{\M}(Spec \, k; \Q(0)) =\Q$ now implies $H^0_{\M}(X; \Q(0)) = \Q(0)$. (Observe also
that the hypothesis $X$ is connected is used only in proving this last condition.) i.e.
We have verified that the DGA $A= A_X=Q_X^{mot}$ associated to the motivic complex of $X$
is connected in the sense of ~\ref{hyp.dga.1}. This proves the first statement. The last statement follows by making use of
the localization sequence in motivic cohomology:
\vskip .3cm
$... \ra H^i_{\M}(Y, {\mathbb Q}(j)) \cong H^{i+2c}_{\M, Y}(X, {\mathbb Q}(j+c)) \ra H^{i+2c}_{\M}(X, {\mathbb Q}(j+c)) \ra H^{i+2c}_{\M}(U, {\mathbb Q}(j+c)) \ra H^{i+1}_{\M}(Y, {\mathbb Q}(j)) \ra \cdots $
\vskip .3cm \noindent
where $c$ is the codimension of $Y$ in $X$. In case $Y$ is a $k$-rational point, clearly its motivic cohomology is
trivial in negative degrees, in view of our hypotheses. Observe that, in general, $Y$ is a projective scheme satisfying the hypotheses above. Therefore, we observe that $H^i_{\M}(Y, {\mathbb Q}(j)) = 0$ for $i<0$ and hence that, for $i+2c<0$, the  map
$  H^{i+2c}_{\M}(X, {\mathbb Q}(j+c)) \ra H^{i+2c}_{\M}(U, {\mathbb Q}(j+c))$ is injective. Moreover
since $c>0$, $i+1<0$ if $i+2c<0$ and therefore, $H^{i+1}_{\M}(Y, {\mathbb Q}(j))=0$ as well. These prove that
$H^{i+2c}_{\M}(U, {\mathbb Q}(j)) =0$ if $i<-2c$. The same observations again show that if
$i=-2c$, the map $H^{i+2c}_{\M}(X, {\mathbb Q}(0)) \ra H^{i+2c}_{\M}(U, {\mathbb Q}(0))$ is an 
isomorphism. The required conclusion now follows.
\end{proof}
\vskip .3cm
The DGA $A$ has a $1$-minimal model, $i:A\langle 1 \rangle \ra A$. (Recall (see \cite[Part IV, section 2]{km}) a connected DGA $B$  over the rationals $\Q$ and provided  with an augmentation $B \ra \Q$ 
is said to be
minimal if it  is a free graded $\Q$-module with decomposable differential $:d(B) \subseteq (I(B))^2$ where $I(B)$ is the augmentation ideal of $B$. $B<1>$ is the sub-DGA of $B$ generated by the
elements of degree $\le 1$ and their differentials. The $1$-minimal model of a DGA $A$ is
a composite map $B \langle 1 \rangle \subseteq B \ra A$ with the last map a quasi-isomorphism and with
$B$ minimal.) The map $i$ induces an isomorphism on
$H^1$ and is injective on $H^2$. We {\it say $A$ is a $K(\pi, 1)$} if $i$ is a quasi-isomorphism.
\begin{theorem} (See \cite [Part IV, Theorem 1.3]{km}.)
The derived category of bounded below chain complexes in ${\mathcal H}_A$ is equivalent to
the derived  category ${\mathbb D}_{-}({A \langle 1 \rangle})$. \end{theorem}
\begin{definition} (The category of  mixed Tate motives over $X$.)
Let ${\chi}^{mot}_X$ denote the Hopf algebra $H^0(\bar B A)$. The category of
(rational)  mixed Tate motives over $X$, denoted ${\mathcal M}{\mathcal T}{\mathcal F}(X)$, will be defined to be the category of finite dimensional co-modules over ${\chi}^{mot}_X$.
\end{definition}
\begin{theorem} 
\label{tatemotives}
 If the DGA $A$ is connected (in the sense of ~\ref{hyp.dga.1}), ${\mathcal M}{\mathcal T}{\mathcal F}(X)$ is equivalent to the category ${\mathcal F}{\mathcal H}_A$. In particular, this holds for the following classes of smooth quasi-projective varieties assuming the Beilinson-Soul\'e conjecture (see above)
holds for the rational motivic cohomology of $Spec \, k$, for example if $k$ is a number field:
\vskip .3cm 
(i) all smooth (connected) projective
linear varieties over $k$
\vskip .3cm
(ii) any of the varieties over $k$ appearing in the list in Examples ~\ref{egs} which are also connected, projective and smooth
\vskip .3cm
(iii) any quasi-projective variety $U$ (over $k$) of the form $X-Y$, where $X$ and $Y$ are smooth projective  varieties both as in (i) or (ii) and $Y$ is closed in $X$.
\end{theorem}
\begin{proof} The proof is clear in view of Theorem ~\ref{Mixed.Tate.1}  and Corollary ~\ref{Mixed.Tate.2}. \end{proof}
\vskip .3cm
Let $\Q(r)$ be the copy of $\Q$ concentrated in bi-degree $(0,r)$ and regarded as a
representation of $\gamma _A$ in the obvious manner.
\begin{corollary} (See \cite[Part IV, Corollary 1.4]{km}.)
If $A=A_X = \Q_X^{mot}$ is a $K(\pi, 1)$, then $Ext^q_{{\mathcal M}{\mathcal T}{\mathcal F}(X)}(\Q, \Q(r)) =
Ext^q_{{\mathcal H}_A}(\Q, \Q(r)) \cong H^q(A(r)) = H^q_{\mathcal M}(X, r) = CH^r(X, 2r-q; \Q)$.
\end{corollary} 


\section{\bf Classical cohomology operations }
{\rm The results of this section follow readily by invoking standard results (see
for example,  \cite{may} or \cite{hsch}) which deduce the existence of cohomological operations
from the existence of an $E_{\infty}$-structure on complexes defining cohomology. However, several nice features of these operations (and hence our constructions) need to be clarified. 
\begin{itemize}
\item The operad
$\{N{\mathbb Z}(E\Sigma _n)|n\}$ is a {\it classical operad} in the sense that the
homology of the complexes $\{N{\mathbb Z}(B\Sigma _n) = N{\mathbb Z}(E\Sigma _n/\Sigma _n)|n\}$ is classical, i.e.
in particular there are not enough classes in the homology of the above complexes to define the motivic operations of Voevodsky. (In fact all the classes have weight $0$.)
\item
However, in \cite{broj} 
we will pursue the relations between these operations and the motivic operations in great detail.  We will see there, that  the motivic operations of Voevodsky and the classical operations considered here differ by multiplication by suitable powers of the Bott element.
\item
 Another interesting feature of our construction is that 
it provides cohomology operations even when $l=p=char(k)$: these are also explored in detail in
\cite{broj}. (The existence of such operations was left open in \cite[section 3, p. 73]{vv2}.)
\item
It also needs to be pointed out that our operations are {\it not} bi-stable, i.e. do not commute 
with weight-suspension, but only with respect to degree-suspension in the sense made precise in the theorem below. 
\end{itemize}}
\begin{remark} In \cite{ep} a purely homological-algebraic technique to defining cohomology operations is considered. However, this requires that the cohomology be with respect to
a {\it sheaf of strictly associative and {\it commutative} algebras}. Therefore, while this approach readily applies to produce cohomology operations in \'etale cohomology with respect to
the sheaves $\{\mu_l|(r)|r\}$, $l \ne char(k)$, it does not apply to motivic cohomology (computed on the Zariski or Nisnevich sites).  Our constructions use nothing more than the existence of an $E_{\infty}$-structure on the
motivic complex.
\end{remark}
\vskip .3cm 
{\rm 
 Let $X$ denote a smooth  separated scheme of finite type over the base field $k$. (We may assume  this is quasi-projective for the sake of simplicity.)
 Let $l$ denote a fixed prime (not necessarily different from the characteristic of $k$) and let ${\underline {\mathbb Z}}_X^{mot}$ 
(${\underline {\mathbb Z}}^{et}_X$) denote the motivic
 complex associated to $X$ restricted to the big Nisnevich site of $X$ (the big \'etale site of $X$, \res.)}
\subsection{}
\label{motivic.1}
 Let ${\mathcal A}_X = {\underline {\mathbb Z}}_X^{mot}{\underset {\mathbb Z} \otimes}{\mathbb Z}/l.{\mathbb
Z}$ and let $A_X = R\Gamma (X, {\mathcal A}_X)$ where the derived
functor is taken on the Zariski (or Nisnevich site). We let
${\mathcal A}_{Spec \, k}$ and $A_{Spec \, k}$ denote the
corresponding objects when $X = Spec \, k$.  
Recall that
${\mathcal A}_X$ and ${\mathcal A}_{Spec \, k}$ are  sheaves of
graded $E_{\infty}$-algebras over the $E_{\infty}$-operad defined in 
the last section. Therefore,
$A_{X}$ and $A_{Spec \, k}$ are  now graded
$E_{\infty}$-differential graded algebras, so that $A_X
= {\underset {r} \oplus} A_X(r)$. Moreover
$H^i_{\mathcal M}(X; {\mathbb Z}/l(r)) = H^i(A_{X}(r))$
which will be isomorphic to  $ CH^{r}(X, 2r-i; {\mathbb Z}/l)$
when $X$ is assumed to be quasi-projective.
\vskip .3cm
\subsection{}
\label{etale.1}
Let ${\mathcal A}_{X,et} = {\underline {\mathbb
Z}}_X^{et}{\underset {\mathbb Z} \otimes}{\mathbb Z}/l.{\mathbb
Z}$ and let $A_{X, et} = R\Gamma (X, {\mathcal A}_{X, et})$ where the derived
functor is taken on the \'etale site.  Now
$H^i_{ et}(X; {\mathbb Z}/l(r)) = H^i(A_{X,et}(r))$.

\subsection{The motivic and \'etale derived categories associated to a scheme}
{\rm For the purposes of 
this section 
we will define this as follows. Recall $({\rm {Sm}}/k)_{Nis}$ ($({\rm {Sm}}/k)_{et}$)
denote the big Nisnevich (\'etale, \res) site of all smooth separated schemes of finite type over the given field $k$. 
For the most part, we will restrict to a fixed smooth scheme $X$. We consider for each smooth separated scheme $X$ of finite type over $k$,
the big Nisnevich site $X_{Nis}$ and the corresponding big \'etale site
$X_{et}$: we may denote either of these generically by $X_{st}$. 
Unless the distinction is important we will continue to
denote both the complexes $\A_X$ and $\A_{X_{et}}$ by $\A_X$ itself.
We consider unbounded (co-chain)
complexes of sheaves $M$ of $Z/l$-vector spaces on the site $X_{st}$. 
We consider the corresponding homotopy category
and the mod$-l$ motivic derived category will be the localization of
this homotopy category by inverting maps that are quasi-isomorphisms.
This category will be denoted $D(X)$. (
We skip these details about the derived category $D(X)$ as they are available in the literature.)
The external hom (internal hom) in this 
category will be denoted $Ext(\quad, \quad )$ ($\RHom(\quad, \quad)$, \res). The internal hom $\RHom(\quad, \quad )$ 
may be made functorial by restricting to cell-$\A$-modules and then by applying the Godement resolution on the second
argument. 
(Observe that $Ext (M, N) =H^0( R\Gamma (X, \RHom(M, N)))$.)
We define the mod$-l$ 
cohomology of an object $M \eps D (X)$ with weight $r$ to be 
$Ext^* (M; \A_X(r))$. This 
will be denoted $H ^*(M; {\mathbb Z}/l(r)) = H^{*,r} (M; {\mathbb Z}/l)$.
\vskip .3cm 
Let $D ^{\le 0}(X)$ denote the full sub-category of $D (X)$ consisting of (co-chain) complexes
$K$ that are trivial in positive degrees. By identifying such co-chain complexes with
chain complexes that are trivial in negative degrees, one may see that that the
derived category $D ^{\le 0}(X)$ is equivalent to the derived category of simplicial
Abelian sheaves on $X_{st} $.  Recall that for any simplicial Abelian sheaf 
$F$ there
is a diagonal map $\Delta: F \ra F \otimes F$; this 
(together with the Alexander-Whitney map and the equivalence between simplicial
Abelian sheaves and complexes of sheaves trivial in negative degrees)
induces a diagonal map
$\Delta : F \ra F \otimes F$, $ F \eps D^{\le 0}(X)$. 
\subsubsection{}
\label{pairing.5.1}
Observe (making use of the above diagonal map) that if $M \eps D^{\le 0}(X)$,
$\RHom (M, \A_X)$ has the obvious induced structure of a sheaf of $E_{\infty}$-algebras
over the operad $\{BE(n)|n \ge 0\}$.  (The required pairings are defined as the
composition 
\vskip .3cm
$BE(n) \otimes \RHom(M, \A_X)^{\overset n \otimes} {\overset
{id \otimes eval} \ra} BE(n) \otimes \RHom(M^{\overset n \otimes}, \A_X^{\overset n \otimes}) {\overset {id \otimes \Delta^*} \ra } BE(n) \otimes 
\RHom(M, \A_X^{\overset n \otimes})$
\vskip .3cm
$ {\overset {\RHom (M, \theta_n)} \ra }
\RHom(M, \A_X)$. 
\vskip .3cm \noindent
The last map is defined by its adjoint: $BE(n) \otimes M \otimes \RHom(M, \A_X^{\overset n \otimes}) {\overset {id \otimes eval} \ra } BE(n) \otimes \A_X^{\overset n \otimes} {\overset {\theta _n} \ra} \A_X$.)
Hence one obtains a graded ring structure on
${\underset r \oplus}H ^*(M; {\mathbb Z}/l(r))$. 
Moreover, if ${\mathbb Z}/l(0)$ denotes
the mod-$l$ motivic complex of weight $0$,  
$\RHom ({\mathbb Z}/l(0), \A_X) \simeq \A_X$
and there is a natural pairing $\RHom (M, \A_X) \otimes \RHom ({\mathbb Z}/l(0), \A_X) \ra \RHom (M, \A_X)$ that is compatible with the above algebra structure on 
\vskip.3cm \noindent
$\RHom (M, \A_X)$. 
\vskip .3cm 
Recall the complex ${\mathbb Z}/l(0)$ is the complex with the constant 
sheaf ${\underline Z}/l$
in degree $0$ and trivial elsewhere in both the motivic and the e\'tale cases.
 Therefore, ${\mathbb Z}/l(0)[i]$ is the
complex concentrated in degree $-i$ where it is the constant sheaf 
${\underline Z}/l$: tensoring with this complex defines the degree-suspension 
$S^i_{deg}$. One may
also obtain the following characterization of the degree-suspension (or the simplicial suspension):
$Ext^* (S^1_{s}M, K) \cong Ext^* (M, K[1])$, for $M, K \eps D (X)$.
We define the Tate suspension (in the motivic case), 
$S^1_{t}M$ by ${\mathbb Z}_{tr}({\mathbb A}^1-0) \otimes M$.
More precisely, we may make use of the pairing
 $\RHom({\mathbb Z}_{tr}({\mathbb A}^1-0), {\mathbb Z}(1)) \otimes 
 \RHom(M, \A_X(r)) \ra \RHom({\mathbb Z}_{tr}({\mathbb A}^1-0) \otimes M, \A_X(r+1))$ to define
 the Tate suspension of the motivic cohomology of $M$. 
The composition of these two suspensions may be effected by tensoring with the
canonical class $T \eps H^2({\mathbb P}^1; {\mathbb Z}/l(1))$. We denote the
composite suspension of $M$ by $S^1_TM$.
 Now we obtain the
natural isomorphisms for any $M \eps D^{\le 0 }(X)$:}
\be \begin{equation}
\begin{split}
\label{suspension.1}
H^{n} (M; {\mathbb Z}/l(r)) \cong H^{n+1}(S^1_{s}M; {\mathbb Z}/l(r)) \quad \\
H^{n} (M; {\mathbb Z}/l(r)) \cong H^{n+1}(S^1_{t}M; {\mathbb Z}/l(r+1)) \quad and \\
H^{n} (M; {\mathbb Z}/l(r)) \cong H^{n+2}(S^1_{T}M; {\mathbb Z}/l(r+1)) \end{split} \end{equation}
\vskip .3cm \noindent
\begin{remarks}
\label{Tate.suspension.etc}
\vskip .3cm
1. Observe that the second isomorphism in the motivic case shows 
$H^1(\RHom({\mathbb Z}_{tr}({\mathbb A}^1-0), {\mathbb Z}(1))) \cong
{ Z}$. Let $\tau$ denote the canonical class corresponding to $1 \eps { Z}$: clearly the 
Tate suspension in the motivic case may be effected by tensoring with this class.
\vskip .3cm
2. All of the above discussion in the motivic case applies equally well when the Nisnevich site is
replaced by the Zariski site.
\end{remarks}
{\rm Throughout the following discussion $H^*$ will denote either motivic or 
\'etale cohomology. We define {\it bi-stable mod-l cohomology 
operations of bi-degree $(i, j)$} 
to be sequences of natural transformations $\{H^n(\quad , {\underline {\mathbb Z}}/l(r)) \ra H^{n+i}(\quad , {\underline {\mathbb Z}}/l(r+j))|n, r\}$
on $D^{\le 0} (X)$ and which are contravariantly functorial in $X \eps ({\rm {Sm}}/k)$.
In view of the suspension-isomorphisms above, these are determined by their restrictions
to $\{H^{2n}(\quad , {\underline {\mathbb Z}}/l(n))|n\}$. 
\vskip .3cm
Recall there are {\it Bockstein} homomorphisms 
$\beta:H^n(M; {\underline {\mathbb Z}}/l(r)) \ra H^{n+1}(M; {\underline {\mathbb Z}}/l(r))$
which are defined in the usual manner as the boundary homomorphism associated to the
short-exact sequence: $0 \ra  {\underline {\mathbb Z}}/l(r) \ra 
{\underline {\mathbb  Z}}/l^2(r) \ra {\underline {\mathbb Z}}/l(r) \ra 0$. These are
clearly bi-stable cohomology operations.
\vskip .3cm
One of the main results in this section is the following theorem, which shows the
existence of classical motivic and \'etale cohomology operations for all primes $l$.}
\vskip .3cm
\begin{theorem}
\label{coh.operations.l.2}
There exist operations $Q^s:H^q(X, {\mathbb Z}/l(t)) \ra 
H^{q+ 2s(l-1)}(X,
{\mathbb Z}/l(l.t))$ and 
\vskip .3cm \noindent
$\beta Q^s:H^q(X, {\mathbb Z}/l(t)) \ra H^{q+ 2s(l-1)+1}(X, {\mathbb
Z}/l(l.t))$. 
\vskip .3cm \noindent 
These operations
satisfy the following properties: \vskip .3cm \noindent 
(i) {\it Contravariant functoriality}: if $f:X \ra Y$ is a map between
smooth separated schemes of finite type over $k$, $f^*\circ Q^s = Q^s \circ f^*$ 
\vskip .3cm \noindent
(ii) Let $x \eps
H^q(X, {\mathbb Z}/l(t))$. $Q^s(x) =0$ if $2s>q$, $\beta Q^s(x) = 0$ if $2s \ge q$ and 
if ($q=2s$), then
 $Q^s(x) = x^l$. 
\vskip .3cm \noindent 
(iii) If $\beta$ is the Bockstein,
$\beta \circ Q^s = \beta Q^s$. \vskip .3cm \noindent 
(iv) {\it Cartan formulae}: {\it For all primes $l$}, $Q^s(x \otimes y) = {\underset
{i+j =s} \Sigma}Q^i(x) \otimes Q^j(y) $ and \vskip .3cm \noindent
$\beta Q^s(x \otimes y) = {\underset {i+j =s} \Sigma } \beta
Q^i(x) \otimes Q^j(y) + Q^i(x) \otimes \beta Q^j(y)$ \vskip .3cm
\vskip .3cm \noindent
(v) {\it Adem relations} For each pair of integers $i \ge 0, j \ge 0$, we let $(i, j) = \frac{(i+j)!}{i!j!}$ with the
convention that $0!=1$. We will also let $(i, j) =0$ if $i<0$ or $j<0$. (See \cite[p. 183]{may}.) With this terminology we obtain:
\vskip.3cm If ($l>2$,  $a<lb$, and $\epsilon =0,1$) or if ($l=2$, $a <lb$ and $\epsilon =0$) one has 
\vskip .3cm 
\be \begin{equation} \beta ^{\epsilon} Q^a  Q^b = \Sigma_i (-1)^{a+i}
(a-li, (l-1)b-a+i-1) \beta ^{\epsilon} Q^{a+b -i} Q^i
\end{equation} \ee
\vskip .3cm \noindent
where $\beta ^{0}Q^s = Q^s$ while $\beta ^1 Q^s = \beta Q^s$. If $l>2$, $a \le lb$ and $\epsilon =0, 1$, one also has
\be \begin{equation}
\begin{split}
\beta ^{\epsilon} Q^a \beta Q^b = (1- \epsilon) \Sigma _i(-1)^{a+i}
(a-li, (l-1)b-a +i-1)\beta Q^{a+b-i} Q^i \\
- \Sigma _i(-1)^{a+i}(a-li-1, (l-1)b-a+i)
\beta ^{\epsilon} Q^{a+b-i} \beta Q^i \end{split} \end{equation} \ee
\vskip .3cm \noindent
(vi) More generally, for any $M\eps D^{\le 0}(X)$, there exist 
cohomology operations 
\vskip .3cm
$Q^s:H^q(M; {\mathbb Z}/l(t)) \ra H^{q+ 2s(l-1)}(M; {\mathbb Z}/l(l.t))$ and 
\vskip .3cm
$\beta Q^s:H^q(M; {\mathbb Z}/l(t)) \ra H^{q+ 2s(l-1)+1}(M; {\mathbb
Z}/l(l.t))$ 
\vskip .3cm \noindent
satisfying the  properties (ii) through (v). Moreover, if $f:M' \ra M$ is a map
in $D^{\le 0}(X)$, the operations $Q^s$ commute with pull-back by $f$. They also commute with the simplicial suspension isomorphism 
in (~\ref{suspension.1}).
\vskip .3cm
(vii) The operation $Q^s$ commutes with change of base fields and also with the
higher cycle map into mod$-l$ \'etale cohomology.
\end{theorem}
\begin{proof}  Let $\pi$ denote the
cyclic group ${Z}/l$ and let $\{e_i|i\}$ form a ${Z}/l$-basis for $H^*(B\pi; Z/l)$. 
We will let $M \eps D^{\le 0}(X)$ and define cohomology operations on
$H^*(M; {\mathbb Z}/l(t)$. For a  smooth scheme $X$, we obtain cohomology operations on 
$H^q(X, {\mathbb Z}/l(t))$ by taking $M = {\mathbf Z}/l_X$, the constant sheaf on $X$ with stalks isomorphic to the integers. Recall that $\A_X$
is an $E_{\infty}$-dga over the $E_{\infty}$-operad $\{BE(n)|n\}$. Since the cohomology operations are assumed to be stable under simplicial suspension as in (ii), it suffices to define these on classes $x \eps H^{2q}(M; {\mathbb Z}/l(t))$. Therefore, one obtains
the existence of cohomology operations $Q^s$ which are defined as follows (see
\cite[p. 161]{may}): if $l=2$, we let
\be \begin{equation}
\label{Q2}
\begin{split}
Q^{s}(x) = \bar\theta_{*}(e_{(2s-2q)} \otimes x^{l}),\\
\beta Q^{s}(x) =  \bar \theta_{*}(e_{(2s-2q)+1} \otimes x^{l})
\end{split}
\end{equation} \ee
\vskip .3cm \noindent 
and if $l>2$, we let: 
\be \begin{equation}
\label{Ql}
\begin{split}
Q^{s}(x) = (-1)^{q-s}  \bar\theta_{*}(e_{(2s-2q)(l-1)} \otimes x^{l}),\\
\beta Q^{s}(x) = (-1)^{q-s} \bar \theta_{*}(e_{(2s-2q)(l-1)+1} \otimes x^{l})
\end{split}
\end{equation} \ee
\vskip .3cm \noindent 
In \cite[p. 161]{may}, an extra sign is introduced. We have gotten rid of
this by including this sign into the choice of the basis elements $\{e_i|i\}$ which form a basis for $H_*(B\pi; { Z}/l)$. 
Observe that $e_{(2s-2q)(l-1)}$ has degree $(2s-2q)(l-1)$ so that the total degree
 of $e_{(2s-2q)(l-1)} \otimes x^l $ is $2q+2s(l-1)$. Since the weight of $x^l$ is
$t.l$ and $\theta_*(e_{(2s-2q)(l-1)} \otimes x^{l})$ leaves the weight unchanged,
$Q^s(x)$ has weight $t.l$.
The map $\bar \theta_*$ is the map 
\be \begin{equation}
\label{bar.theta}
\bar \theta_*:H_*(B\pi; { Z}/l) \otimes H^*(M; {\mathbb Z}/l(r))^{\powerl} 
\ra  H^*(M; {\mathbb Z}/l(r)), \quad \bar \theta_* (\bar e \otimes \bar x^{\powerl}) = [\theta (e \otimes x^{\powerl})]
\end{equation} \ee
\vskip .3cm \noindent
where  $e$ ($x$) denotes a cycle
representing the cohomology class $\bar e$ ($\bar x$, \res) and $[ z ]$ denotes
the cohomology class represented by the cycle $z$.
\vskip .3cm
Now all assertions except for the 
last assertion are immediate consequences of standard 
results on  cohomology operations on algebras over ${\rm E}_{\infty}$-operads: see \cite{may} or \cite{hsch}. (Recall that $e_i$ is defined only for $i\le 0$; therefore we let $e_i=0$ for $i>0$ so that $Q^s(x) = 0 $ if $q<2s$.)
The action of the  operad $\{BE(p)|p\}$
on the complex $\A_X$ was shown to be functorial in the base field $k$; therefore
the operation $Q^s$ commutes with respect to change of base fields. Moreover,
as observed in Remark ~\ref{Tate.suspension.etc}, the action of the operads
$\{BE(n)|n \ge 0\}$ is compatible on the complexes $\A_X$ and
$\A_{X_{et}}$; therefore, the operation $Q^s$ is also compatible with respect to
the higher cycle map from mod-$l$ motivic to mod-$l$ \'etale cohomology.. 
\end{proof}
\begin{remark} The above operations cannot commute
with the Tate suspension as one may see from elementary weight considerations. 
Therefore they are {\it not} bi-stable cohomology operations.
Moreover, observe that  $Q^0$ is  not necessarily the identity and therefore 
$\beta Q^0$ is not necessarily $\beta$.  Observe also that some of the
above cohomology operations may be trivial owing to the fact that the weight may be
high enough. (Recall: $H^i_{\M}(X; {\mathbb Z}/l(j)) \cong CH^{j}(X, {\mathbb Z}/l; 2j-i) =0$
if $j>dim (X) +2j-i$.) This shows that the motivic cohomology operations of Voevodsky {\it cannot} be deduced from the classical operations considered above. 
However, the classical operations on \'etale cohomology indeed can be deduced from the motivic cohomology operations : this is discussed in detail in the 
 paper \cite{broj}.
\end{remark}
{\rm We will denote the operations $Q^s$ considered above on \'etale cohomology with {\it $l$ different the characteristic $p$} as $Q^s_{et}$.
When the base field is separably closed, one may identify $\mu_l(r)$ with the constant
sheaf ${ Z}/l$; in this case, therefore, the weights are irrelevant, and we 
obtain cohomology operations in the usual sense, once it is shown that $Q^0_{et} =id$. 
(This is proved below.)
For example, if the base field is
the field of complex numbers, the operations we obtain 
identify with the usual cohomology operations in mod-$l$
singular cohomology (once the latter is identified with mod$-l$ \'etale cohomology).}
\begin{proposition} If the base field is separably closed, 
 the operation $Q^0_{et} =id$.
\end{proposition}
\begin{proof} Since the Tate suspension is irrelevant now, 
$Q^0$ commutes with the simplicial suspension and is contravariantly
functorial we may reduce to checking this when $M= $ the constant sheaf ${ Z}/l$.
In this case, $Q^0_{et}(\alpha) = \alpha ^l = \alpha $, for any $\alpha \eps { Z}/l$.
This proves the proposition. \end{proof}
\vskip .5cm

\section{\bf Appendix: chain complexes and operads}
\label{DP}
{\rm  Let $\mathbf A$ denote an Abelian category; a chain complex $K$ in $\mathbf A$ will denote a sequence $K_i \eps {\mathbf A}$ provided with
maps $d : K_i \ra K_{i-1}$ so that $d^2 = 0$. Let $\oC_0({\mathbf A})$ denote the
category of chain complexes in ${\mathbf A}$ that are trivial in negative degrees. One defines
the de-normalizing functor:
\be \begin{equation}
\label{DN}
DN: \oC_0({\mathbf A}) \ra (simplicial\ objects \ in \ {\mathbf A})]
\end{equation} \ee
as in \cite[pp.8-9]{ill} so that 
$DN$ will be inverse to  the
functor 
\be \begin{equation}
 \label{N}
N: (simplicial \ objects \ in \ {\mathbf A}) \ra \oC_0({\mathbf A})
\end{equation} \ee
 defined by $(NK)_n = {\underset {i \neq 0} \cap } ker (d_i: K_n \ra K_{n-1})$ with $\delta : (NK)_n \ra (NK)_{n-1}$ induced by $d_0$. We will also often consider the functor $C: (simplicial \ objects \ in \ {\mathbf A}) \ra \oC_0({\mathbf A})$ defined by $K \mapsto $ the chain complex which in degree $n$ is $K_n$ and where the 
differential $d: C(K)_n \ra C(K)_{n-1}$ is given by $d = \Sigma _i (-1)^i d_i$. Given a chain complex $K$, trivial  in negative degrees with differentials of
degree $-1$, one may form the corresponding (co-)chain complex trivial in positive degrees with differentials of degree $+1$ by re-indexing $K$ in the obvious manner. This functor composed with the functor $C$ above will be denoted $C^{\bullet}: (simplicial \ objects \ in \ {\mathbf A}) \ra \oC({\mathbf A})$ where
the last is the category of (co-)chain complexes in the Abelian category ${\mathbf A}$.
\vskip .3cm
Given two positive integers $p$ and $q$, a $(q, p)$-{\rm {shuffle}} $\pi$ is a permutation of $(1,..., p+q)$
so that $\pi(i)< \pi(j)$ for $1\le i<j \le q$ and for $q+1 \le i <j \le q+p$. We let $\mu=$ the
restriction of $\pi$ to $ (1,..., q)$ and $\nu( j)=\pi(j+q)$, $1 \le j \le p$. Clearly $ \pi$ is
determined by $(\mu, \nu)$ and therefore, we will identify $\pi$ with the pair $(\mu, \nu)$. The set of all $(q, p)$-shuffles is in one-one correspondence with the
set of all strictly increasing maps $(\phi, \psi):[p+q] \ra [q] \times [p]$: the
correspondence is given by sending a shuffle $(\mu, \nu)$ to $(\phi, \psi)$ where
$\phi(x) = $ the cardinality of the set $\{1 \le i \le q| (\mu, \nu)(i) \le x\}$ and $\psi(x) = $ the cardinality of the set $\{q+1 \le i \le q+p| (\mu, \nu)(i) \le x\}$. Each such map $(\phi, \psi)$ (and hence each shuffle map) defines an isomorphism of schemes
\be \begin{equation}
\label{shuffle.iso}
\Delta[p+q] \ra \Delta[q] \times \Delta [p]
\end{equation} \ee
\vskip .3cm \noindent
by the formula: $(t_0, \cdots, t_{p+q}) \mapsto t_0(v_{\phi(0)}, v_{\psi(0)}) + \cdots +t_{p+q}(v_{\phi(p+q)}, v_{\psi(p+q)})$.
\vskip .3cm
\subsubsection{}
\label{shuffle.Esigma}
Let $\Sigma _k$ denote the symmetric group on $k$-letters and let $E\Sigma _k$ denote the simplicial group defined
by the bar construction. Observe that a $p$-simplex of $E\Sigma _k$ is given by a sequence $(\sigma _0, \cdots, \sigma _p)$, $\sigma _i \eps \Sigma _k$. Let $s_p = (\sigma _0, \cdots, \sigma _p)$ and $s_q = (\tau _0, \cdots, \tau _q)$
denote a $p$-simplex and a $q$-simplex of $E\Sigma _k$. 
Then, viewing $s_p$ as a map $\Delta [p] \ra E\Sigma_k$, $s_q$ as a map $\Delta [q] \ra E\Sigma _k$ of simplicial sets, $\phi$ as an element of $\Delta[q]_{p+q}$ and $\psi$ as an element of $\Delta[p]_{p+q}$, each  map $(\phi, \psi)$ (as above) associates a $p+q$-simplex $(s_q(\phi(0)) \circ s_p(\psi(0)), \cdots, s_q(\phi (p+q)) \circ s_p(\psi(p+q)))$ of $E\Sigma _k$ to the  $(q,p)$-simplex $(s_q,  s_p)$. Here $\circ$ denotes composition in $\Sigma _k$. We denote this product by
\be \begin{equation}
\label{shuffle.product}
(\phi, \psi)^*(s_q \times s_p)
\end{equation} \ee
\vskip .3cm
Given two simplicial objects $K$ and $K'$, one obtains a map called the shuffle map (where $C$ denotes the functor  considered above)
\be \begin{equation}
s={\rm {shuffle}}: K_q \otimes K'_p = C(K)_q \otimes C(K')_p \ra C(K \otimes K')_{p+q} 
\end{equation} \ee
which is defined by
\be \begin{equation}
{\rm {shuffle}}(k_q \otimes k'_p) = \Sigma _{(\mu, \nu)} (-1)^{\sigma (\nu)}(s_{\mu _p} \circ ... \circ s_{\mu _1}(k_q)  \otimes s_{\nu_q} \circ ... \circ s_{\nu _1}(k_p'))
\end{equation} \ee
where the sum is over all $(q, p)$-{\rm {shuffle}}s $(\mu, \nu)$ and where $\sigma(\nu)$ is the signature of the permutation $\nu$. }
\begin{proposition}
\label{shuffle.props}
 The functor $C$ is compatible with pairings. Moreover, the
shuffle map above is strictly associative and strictly commutative.
\end{proposition}
\subsection{Operads}
{\rm Basic definitions of operads and algebras over operads  may be found in \cite[Part I, section 1]{km} or \cite{hsch}. The same definitions apply with minor modifications to the unital symmetric monoidal category
$\oC(\S, \R)$ as in section 2. We summarize some of these definitions very briefly here, mainly for the sake of the reader unfamiliar with operads.  
\vskip .3cm
An {\it associative operad} (or simply operad)
${\mathcal O}$ in $\oC(\S, \R)$  is given by a sequence $\{{\mathcal O}(k)|k \ge 0\}$ of objects in
$\oC(\S, \R)$ along with the following data: for every integer $k \ge 1$ and every sequence $(j_1,..., j_k)$ of non-negative integers so that $\Sigma _{l=1}^{l=k} j_l = j$ there is given a map
$\gamma _k:{\mathcal O}(k) \otimes {\mathcal O}(j_1) \otimes ... {\mathcal O}(j_k) \ra {\mathcal O}(j)$ so that certain associativity diagrams commute. In addition one is provided with a unit map $\eta : u \ra {\mathcal O}(1)$ so that it acts as a unit for the compositions $\gamma: {\mathcal O}(1) \otimes {\mathcal O}(j) \ra {\mathcal O}(j)$ and for $\gamma: {\mathcal O}(k) \otimes {\mathcal O}(1)^{\otimes \, k} \ra {\mathcal O}(k)$. 
\vskip .3cm
A {\it commutative operad} is an operad as above provided with an action
by the symmetric group $\Sigma _k $ on each ${\mathcal O}(k)$ so that the
maps $\gamma$ are compatible with the $\Sigma_k$-action, which is expressed by saying that certain diagrams
as in \cite[section 1]{hsch} (or \cite[Part I, section 1]{km}) commute.
\vskip .3cm
 An operad is an $A_{\infty}$-operad (or {\it acyclic} operad) if
{\it each ${\mathcal O}(k)$ is {\it acyclic}}.
It is an ${\it E_{\infty}}$-operad, if in addition, it is commutative and {\it the given action of $\Sigma_k$ on ${\mathcal O}(k)$ is {\it free}}.
\begin{examples} 
\label{EZ.operad}
1. {\it The classical Eilenberg-Zilber operad}. We will recall the definition of this operad briefly starting with endomorphism operads.
 Consider the functor $\Delta \ra \oC(\S, R)$, defined by $n \mapsto C_*(\Delta [n]_{ss}, R)$, where $\Delta [n]_{ss}$ denotes the 
simplicial set $\{Hom_{\oDelta}([k], [n])| [k] \eps \oDelta \}$. This complex re-indexed so that it is trivial in positive degrees and where 
the differentials are of degree $+1$ will be denoted ${\mathcal Z}$. (The subscript $_{ss}$ is used in $\Delta[n]_{ss}$ to distinguish this
object from the scheme $\Delta [n]$ considered elsewhere in the paper.) 
\vskip .3cm
 We have the  diagonal power functor $d \mapsto {\mathcal Z} ^{\otimes \, d}$ which is the $d$-fold tensor product of the functor ${\mathcal Z}$, i.e.
it sends $[n] \eps \oDelta$ to $C^*(\Delta[n])$.
By convention, the $0$-fold power of ${\mathcal Z}$ is the constant functor at $R$. We define the {\it endomorphism operad} $\End_{{\mathcal Z}}$ of the functor $C^*$ by 
letting 
\be \begin{equation}
\label{end.op.gen}
\End_{{\mathcal Z}}(n) = \Hom_{{\oDelta}}({\mathcal Z}, {\mathcal Z}^ { \otimes \,  n })
\end{equation} \ee
\vskip .3cm \noindent
with $\End_{{\mathcal Z}}(0) = R$. The structure morphisms are defined as follows.  $\theta _{n}:{\mathcal E}nd_{{\mathcal Z}}(n) \otimes {\mathcal E}nd_{{\mathcal Z}}(k_1) \otimes \cdots \otimes {\mathcal E}nd_{{\mathcal Z}}(k_n) \ra {\mathcal E}nd_{{\mathcal Z}}(\Sigma _i k_i)$  
is defined as the composition of $\Hom({\mathcal Z}, {\mathcal Z}^{\otimes \, n}) \otimes \Hom({\mathcal Z}, {\mathcal Z}^{\otimes \, {k_1}}) \otimes \cdots \otimes \Hom({\mathcal Z}, {\mathcal Z}^{\otimes \, {k_n}}) \ra \Hom({\mathcal Z}, {\mathcal Z}^{\otimes ^n}) \otimes \Hom({\mathcal Z}^{\otimes \, n}, {\mathcal Z}^{\otimes^{ \Sigma _i k_i}}) \ra\Hom({\mathcal Z}, {\mathcal Z}^{\otimes \, {\Sigma _i k_i}})$ where both maps are the obvious ones.
 This operad is known to be an acyclic operad: see \cite{hsch} for more details.
\vskip .3cm
2. {\it The trivial operad}. This is obtained by taking each ${\mathcal O}(n) = \R$ provided with the trivial action of $\Sigma_n$.
\end{examples}
\vskip .3cm
\begin{definition}
\label{algebra}
 A differential graded algebra  (or algebra)  $A$ over an operad $\{{\mathcal O}(n)|n \}$ is an object in $\oC(\S, \R)$ provided with
maps $\theta : {\mathcal O}(j) \otimes {A }^ {\otimes \, j} \ra A$ for all $j \ge 0$ that are associative and unital in 
the
sense that the following diagrams commute:
\vskip .3cm
\be \begin{multline}
\label{operad.4}
{\diagram
{\mathcal O}(k) \otimes {\mathcal O}(j_1) \otimes ... \otimes {\mathcal O}(j_k) \otimes {A }^{ \otimes \, j} \stackrel{\gamma \otimes id}{\longrightarrow} \ddto_{{\rm {shuffle}}} &{\mathcal O}(j) \otimes {A }^{ \otimes \, j} \dto_{\theta}\\
& A\\
{\mathcal O}(k) \otimes {\mathcal O}(j_1) \otimes A^{ \otimes \, {j_1}}\otimes ... \otimes {\mathcal O}(j_k) \otimes {A }^{  \otimes \, {j_k}} \overset { id \otimes \theta ^k}{\longrightarrow}  &{\mathcal O}(k) \otimes {A }^{\otimes \, k} \uto_{\theta}
\enddiagram}
\end{multline} \ee
 and
\hskip .5cm
\be \begin{multline}
\label{operad.5}
{\diagram
u \otimes A \rto_{\simeq} \dto _{\eta \otimes id}& A\\
{\mathcal O}(1) \otimes A \urto_{\theta}\enddiagram} \end{multline} \ee
\vskip .3cm \noindent
See \cite[section 1]{hsch} or \cite[Part I, section 1]{km}. We will often refer to such algebras as dgas or dg-algebras over the operad 
$\{{\mathcal O}(n)|n\}$.
If the operad is ${\it A_{\infty}}$, we will refer to the algebra $A$ as an ${\it A_{\infty}}$-algebra, i.e. an $A_{\infty}$-dga.
\end{definition}
\vskip .3cm
A commutative algebra over a commutative operad ${\mathcal O}$ is an $A_{\infty}$ algebra over the
operad ${\mathcal O}$ so that the following diagrams commute:
\vskip .3cm
 \be \begin{multline}
\label{operad.algebra}
{\diagram
{\mathcal O}(j) \otimes A ^{ \otimes \, j} \rrto_{\sigma \otimes \sigma ^{-1}} \drto_{\theta}&& {\mathcal O}(j) \otimes A^{ \otimes \, j } \dlto^{\theta}\\
& A\enddiagram} \end{multline}
 \ee
\vskip .3cm \noindent
If, in addition, the operad is $E_{\infty}$, we will refer to the algebra $A$ as an ${\it E_{\infty}}$-algebra, i.e. an $E_{\infty}$-dga. A {\it commutative differential graded algebra}
 or {\it commutative dga} is a commutative algebra over the trivial operad discussed in Example ~\ref{EZ.operad} 2.
\vskip .3cm
{\it A cosimplicial algebra} $A^{\bullet}$ over a commutative operad $\{{\mathcal O}(n)|n\}$ is a functor $\Delta \ra ($algebras over the operad $\{{\mathcal O}(n)|n\}$), i.e. a cosimplicial object in the category of algebras over the
operad $\{{\mathcal O}(n)|n\}$. Given a cosimplicial algebra $A^{\bullet}$, one may take its normalization $N(A^{\bullet})$ to be the total complex of the
double complex obtained by first normalizing $A^{\bullet}$ in the cosimplicial direction. (The total complex of a double complex 
$K=\{K^{n, m}|n, m\}$ will denote the complex $Tot(K) $ defined by $Tot(K)^n = {\underset {n=i+j} \Pi} K^{i, j}$ with the induced boundary map.)}
\begin{proposition} 
\label{norm.operad}
If $A^{\bullet}$ is a cosimplicial algebra over the commutative operad $\{{\mathcal O}(n)|n\}$, its normalization $N(A^{\bullet}) $ is 
an algebra over the operad $\{{\mathcal O}(n)|n\} \otimes \End_{{\mathcal Z}}$.
\end{proposition}
\begin{proof} This is a straightforward extension of \cite[(2.3) Theorem]{hsch}. One first shows that the double complex obtained by 
normalizing $A^{\bullet}$ in the cosimplicial direction is an algebra over the operad 
$\{\{{\mathcal O}^i(n) \otimes \End_{{\mathcal Z}}^j(n)|i, j\}| n\}$ in the category of
double complexes. (Here the superscripts $i$ and $j$ are the degrees of the double complex.) Now one takes the associated double complexes. 
(See \cite[proof of (2.3) Theorem]{hsch} for further details.) The tensor product of operads here is defined simply by taking the
tensor product of the corresponding complexes. In order to see this defines an operad one may consult \cite[4.1]{Moerd}. \end{proof}
\vskip .3cm \noindent
 
\begin{definition} Recall we have assumed that the site $\S$ has {\it enough points}. Therefore one may define {\it the Godement resolution} 
$\{G^nK|n\}$, as a cosimplicial resolution of 
any object $K \eps \oC(\S, \R)$: see \cite{josh2},  for example.
\label{coh.des.MV}
(i) We say  $K \eps \oC(\S, \R)$  has cohomological descent on the
site $\S$, if for each $V$ in $\S$, the natural map $\Gamma (V, K) \ra
\holimD \{ \Gamma (V, {G^nK})|n\}$ is a weak-equivalence.
\vskip .3cm
(ii) Let $\S$ denote the Zariski site of a Noetherian scheme $X$. We say that
a $K \eps \oC(X_{Zar})$  has the Mayer-Vietoris property on $X$, if for any two 
Zariski open sub-schemes $U$, $V$ of $X$, the diagram
\vskip .3cm
$\Gamma (U \cup V, K) \ra \Gamma (U, K) \oplus \Gamma (V, K) \ra \Gamma (U \cap V, K) \ra  \Gamma(U \cup V, K[1])$
\vskip .3cm \noindent
is a distinguished triangle of complexes of abelian groups.
\end{definition}
\vskip .3cm \noindent
The following is a well-known result: see for example, \cite{T}.
\begin{proposition}
\label{coh.des}
 Suppose $X$ is a Noetherian scheme and $K \eps \oC(X_{Zar})$
 having the Mayer-Vietoris property. Then
$K$ has cohomological descent on the Zariski site of $X$. \end{proposition}
\vskip .3cm
In general, for  a site $\S$ with enough points and $K \eps \oC(\S, \R)$ we let $R\Gamma (V, K) = \holimD \{\Gamma (V, G^nK)|n\}$. 
\vskip .3cm


\end{document}